\documentclass[aap,MSNbibl,dvips]{arximspdf}
\usepackage{algorithmic}
\usepackage{mathbh}
\usepackage{ifthen}
\usepackage{xargs}

\doi{10.1214/13-AAP962} 
\volume{24}
\issue{5}
\pubyear{2014}
\firstpage{1767}
\lastpage{1802}

\makeatletter
\newcommand{\widebar}{\overline}

\newcommand{\rright}{\right}
\newcommand{\lleft}{\left}
\newcommand{\rrVert}{\Vert}
\newcommand{\rrvert}{\vert}
\newcommand{\llVert}{\Vert}
\newcommand{\llvert}{\vert}
\newcommand{\varotimes}{\otimes}
\newcommand{\eqref}[1]{(\ref{#1})}
\newcommand{\iint}{\int\!\!\int}
\newcommand{\idotsint}{\int\cdots\int}
\newtheorem{theorem}{Theorem}
\newtheorem{proposition}[theorem]{Proposition}
\newtheorem{lemma}[theorem]{Lemma}
\newtheorem{corollary}[theorem]{Corollary}
\newproclaim{definition}[theorem]{Definition}
\newproclaim{remark}[theorem]{Remark}

\newcommand{\coint}[1]{[#1[}

\newcommand{\ooint}[1]{]#1[}
\newcommand{\ccint}[1]{[#1]}
\newcommand{\alg}[1]{\mathcal{#1}}
\newcommand{\algprod}{\varotimes}
\newcommand{\bmf}[1]{\mathcal{F}(#1)}
\newcommand{\borel}[1]{\mathcal{B}(#1)}
\newcommand{\chunk}[3]{{#1}_{#2}^{#3}}
\newcommand\DDelta[3]{\Delta_{#1}\langle #2\rangle(#3)}
\newcommand\DDeltaantri[3]{\Delta_{#1}\bigl\langle #2\bigr\rangle\bigl(#3\bigr)}

\newcommand{\dlim}{\stackrel{\mathcal D}{\longrightarrow}}
\newcommand{\epart}[2]{\ensuremath{\xi_{#1}^{#2}}}
\newcommand{\eqdef}{\triangleq}
\newcommand{\espcan}{\mathbb{E}_\chi}
\newcommand{\fm}[1]{\mu_{#1}}

\newcommand{\G}{\mathbf{G}}
\newcommand{\ind}[2]{I_{#1}^{#2}}
\newcommand{\K}{\mathbf{K}}
\newcommand{\Lp}[1]{\mathsf{L}^{#1}}
\newcommand{\M}{\mathbf{M}}

\newcommand{\mcp}[2]{\mathcal{M}{\langle #1 \rangle}(#2)}
\newcommand{\mcpantri}[2]{\mathcal{M}{\bigl\langle #1 \bigr\rangle}(#2)}
\newcommand{\mdr}{\mathcal{M}(\set{D},r)}
\newcommand{\meas}[1]{\mathcal{M}(#1)}
\newcommand{\midbar}{\vert}

\newcommand{\modulotxt}[1]{[#1]_r}
\newcommand{\nset}{\mathbb{N}}
\newcommand{\nsetpos}{\mathbb{N}^\ast}
\newcommand{\1}{\mathbh{1}}
\newcommand{\probdoeblin}[2]{\mu_{#1}\langle #2 \rangle}
\newcommand{\pmf}[1]{\mathcal{F}_+(#1)}

\newcommand{\probcan}{\mathbb{P}_\chi}
\newcommand{\probmeas}[1]{\mathcal{P}(#1)}
\newcommand{\Q}{\mathbf{Q}}

\newcommand{\rset}{\mathbb{R}}
\newcommand{\rsetpos}{\mathbb{R}^\ast_+}
\newcommand{\set}[1]{{\mathsf{#1}}}
\newcommand{\transp}[1]{{^t#1}}
\newcommand{\T}{\mathbf{T}}
\newcommand{\upstep}[1]{\Psi\langle #1 \rangle}

\newcommand{\wgt}[2]{\omega_{#1}^{#2}}
\newcommand{\wgtsum}[1]{\Omega_{#1}^N}
\newcommand{\Xcan}{X}
\newcommand{\Ycan}{Y}
\newcommand{\Yfilt}[1]{\mathcal{F}^Y_{#1}}
\newcommand{\yref}{\nu}
\newcommand{\zset}{\mathbb{Z}}

\newcommandx\A[2][1=]{
\ifthenelse{\equal{#1}{}}
{(\textbf{A\ref{#2}})}
{(\textbf{A\ref{#1}--\ref{#2}})}
}
\newcommandx\dens[3][1=,3=]%
{
\ifthenelse{\equal{#1}{}}
{\ifthenelse{\equal{#3}{}}{p(#2)}{p_{#3}(#2)}}
{\ifthenelse{\equal{#3}{}}{p_{#1}(#2)}{p_{#1}^{#3}(#2)}}
}
\newcommandx\esp[3][1=,3=]{{\mathbb E}_{#1}^{#3} (#2)}
\newcommandx\espbig[3][1=,3=]{{\mathbb E}_{#1}^{#3} \bigl(#2\bigr)}
\newcommandx\espBig[3][1=,3=]{{\mathbb E}_{#1}^{#3} \Bigl(#2\Bigr)}
\newcommandx\espBigg[3][1=,3=]{{\mathbb E}_{#1}^{#3} \Biggl(#2\Biggr)}
\newcommandx\cesp[4][1=,4=]{{\mathbb E}_{#1}^{#4} (#2| #3)}
\newcommandx\cespbig[4][1=,4=]{{\mathbb E}_{#1}^{#4} \bigl(#2| #3\bigr)}

\newcommandx\espTxt[3][1=,3=]{
\ifthenelse
{\equal{#3}{}}
{{\mathbb E}_{#1} (#2)}
{{\mathbb E}_{#1} (#2|#3)}
}
\newcommandx\feyn[2][1=]
{
\ifthenelse{\equal{#1}{}}
{\ensuremath{\Phi\langle #2 \rangle}}
{\ensuremath{\Phi\langle #2 \rangle(#1)}}
}

\newcommandx\feynantri[2][1=]
{
\ifthenelse{\equal{#1}{}}
{\ensuremath{\Phi\bigl\langle #2 \bigr\rangle}}
{\ensuremath{\Phi\bigl\langle #2 \bigr\rangle\bigl(#1\bigr)}}
}

\newcommand{\field}[1]{\ensuremath{\mathcal{#1}}}

\newcommandx\filt[3][1=,3=]
{
\ifthenelse
{\equal{#1}{}}
{\ensuremath{{\bar{\phi}}\langle #2 \rangle}}
{\ensuremath{{\bar{\phi}}_{#1} \langle #2 \rangle}}%
}

\newcommandx\filtantri[3][1=,3=]
{
\ifthenelse
{\equal{#1}{}}
{\ensuremath{{\bar{\phi}}\bigl\langle #2 \bigr\rangle}}
{\ensuremath{{\bar{\phi}}_{#1} \bigl\langle #2 \bigr\rangle}}%
}

\newcommandx\pred[3][1=,3=]%
{%
\ifthenelse
{\equal{#1}{}}
{\ensuremath{{\phi}\langle #2 \rangle}}
{\ensuremath{{\phi}_{#1} \langle #2 \rangle}}%
}

\newcommandx\predantri[3][1=,3=]%
{%
\ifthenelse
{\equal{#1}{}}
{\ensuremath{{\phi}\bigl\langle #2 \bigr\rangle}}
{\ensuremath{{\phi}_{#1} \bigl\langle #2 \bigr\rangle}}%
}

\newcommandx\filtvariance[3][1=]{\bar{\sigma}_\chi^{#1} \langle #2\rangle(#3)}
\newcommandx\filtvarianceantri[3][1=]{\bar{\sigma}_\chi^{#1} \bigl\langle #2\bigr\rangle(#3)}

\newcommandx\LGSS[2][1=]{
\ifthenelse{\equal{#1}{}}
{(\textbf{LGSS\ref{#2}})}
{(\textbf{LGSS\ref{#1}--\ref{#2}})}
}
\newcommand{\lleb}{\lambda^{\mathrm{Leb}}}

\newcommand{\Lkh}[1]{\mathbf{L}\langle{#1}\rangle}
\newcommand{\Lkhantri}[1]{\mathbf{L}\bigl\langle{#1}\bigr\rangle}
\newcommandx\llhd[2][2=]{
\ifthenelse{\equal{#2}{}}
{g \langle #1 \rangle}
{g \langle #1 \rangle(#2)}
}
\newcommandx\llhdantri[2][2=]{
\ifthenelse{\equal{#2}{}}
{g \bigl\langle #1 \bigr\rangle}
{g \bigl\langle #1 \bigr\rangle\bigl(#2\bigr)}
}

\newcommandx\pot[2][1=]{
\ifthenelse{\equal{#1}{}}
{G_{#2}}
{G_{#2}(#1)}
}

\newcommandx\partfilt[3][1=,3=]
{
\ifthenelse
{\equal{#1}{}}
{\ensuremath{\bar{\phi}^N \langle #2 \rangle}}
{\ensuremath{{\bar{\phi}}_{#1}^N \langle #2 \rangle}}%
}

\newcommandx\partfiltantri[3][1=,3=]
{
\ifthenelse
{\equal{#1}{}}
{\ensuremath{\bar{\phi}^N \bigl\langle #2 \bigr\rangle}}
{\ensuremath{{\bar{\phi}}_{#1}^N \bigl\langle #2 \bigr\rangle}}%
}

\newcommandx\partpred[3][1=,3=]%
{%
\ifthenelse
{\equal{#1}{}}
{\ensuremath{{\phi}^N\langle #2 \rangle}}
{\ensuremath{{\phi}^{N}_{#1} \langle #2 \rangle}}%
}

\newcommandx\partpredantri[3][1=,3=]%
{%
\ifthenelse
{\equal{#1}{}}
{\ensuremath{{\phi}^N\bigl\langle #2 \bigr\rangle}}
{\ensuremath{{\phi}^{N}_{#1} \bigl\langle #2 \bigr\rangle}}%
}

\newcommand{\Pblock}[2][]
{\ifthenelse{\equal{#1}{}}{\mathbf{L}\langle#2\rangle}{\mathbf{L}^{#1}\langle#2\rangle}
}
\newcommand{\Pblockantri}[2][]
{\ifthenelse{\equal{#1}{}}{\mathbf{L}\bigl\langle#2\bigr\rangle}{\mathbf{L}^{#1}\bigl\langle#2\bigr\rangle}
}
\newcommand{\pblock}[2][]
{\ifthenelse{\equal{#1}{}}{\bolds{\ell}\langle#2\rangle}{\bolds{\ell}^{#1}\langle #2\rangle}
}
\newcommand{\pblockantri}[2][]
{\ifthenelse{\equal{#1}{}}{\bolds{\ell}\bigl\langle#2\bigr\rangle}{\bolds{\ell}^{#1}\bigl\langle #2\bigr\rangle}
}

\newcommandx\likeli[3][1=]{\pi_{#1}\langle #2 \rangle(#3)}
\newcommandx\likeliantri[3][1=]{\pi_{#1}\bigl\langle #2 \bigr\rangle(#3)}

\def\probalone{\ensuremath{\mathbb{P}}}
\newcommandx
\prob[3][1=,3=]{
\ifthenelse
{\equal{#3}{}}
{{\mathbb P}_{#1}(#2)}
{{\mathbb P}_{#1} (#2| #3)}
}

\newcommandx\probBig[3][1=,3=]{
\ifthenelse
{\equal{#3}{}}
{{\mathbb P}_{#1}\Bigl(#2\Bigr)}
{{\mathbb P}_{#1} \Bigl(#2| #3\Bigbr)}
}
\newcommandx\probantri[3][1=,3=]{
\ifthenelse
{\equal{#3}{}}
{{\mathbb P}_{#1}\bigl(#2\bigr)}
{{\mathbb P}_{#1} \bigl(#2| #3\bigr)}
}

\newcommandx\probTxt[3][1=,3=]{
\ifthenelse
{\equal{#3}{}}
{{\mathbb P}_{#1} (#2)}
{{\mathbb P}_{#1} ( #2 | #3)}
}

\newcommandx\secfact[2][1=]{\Lambda_{#1} \langle #2 \rangle(h)}
\newcommandx\secfactub[2][1=]{\Upsilon_{#1} \langle #2 \rangle(h)}
\newcommandx\supnorm[2][1=]{\ensuremath{\llVert #2\rrVert ^{#1}_{\infty}}}
\newcommandx\supnormantri[2][1=]{\ensuremath{\bigl\llVert #2\bigr\rrVert ^{#1}_{\infty}}}
\newcommand{\supnormTxt}[1]{\ensuremath{\|#1\|_{\infty}}}
\newcommand{\supnormTxtantri}[1]{\ensuremath{\bigl\|#1\bigr\|_{\infty}}}
\newcommandx\sequence[3][2=n,3=\zset]{\ensuremath{(#1_{#2})_{#2 \in #3}}}
\newcommandx\variance[3][1=]{\sigma_\chi^{#1} \langle #2 \rangle(#3)}
\newcommandx\varianceantri[3][1=]{\sigma_\chi^{#1} \bigl\langle #2 \bigr\rangle(#3)}

\newcommand{\URoot}[1][]
{\ifthenelse{\equal{#1}{}}{\ensuremath{R}}{\ensuremath{R}_{#1}}}
\newcommand{\VRoot}[1][]
{\ifthenelse{\equal{#1}{}}{\ensuremath{S}}{\ensuremath{S}_{#1}}}
\newcommand{\UCov}[1][]%
{%
\ifthenelse{\equal{#1}{}}{\URoot{^t\URoot}}{\URoot[#1] {^t\URoot[#1]}}%
}
\newcommand{\VCov}[1][]%
{%
\ifthenelse{\equal{#1}{}}{\VRoot{^t\VRoot}}{\VRoot[#1] {^t\VRoot[#1]}}%
}
\makeatother

\begin{document}
\begin{frontmatter}

\title{Long-term stability of sequential Monte Carlo methods under
verifiable conditions}
\runtitle{Long-term stability of sequential Monte Carlo methods}

\begin{aug}
\author[A]{\fnms{Randal} \snm{Douc}\thanksref{both}\ead[label=e1]{randal.douc@it-sudparis.eu}},
\author[B]{\fnms{Eric} \snm{Moulines}\thanksref{both}\ead[label=e2]{eric.moulines@telecom-paristech.fr}}
\and
\author[C]{\fnms{Jimmy} \snm{Olsson}\corref{}\thanksref{jo}\ead[label=e3]{jimmyol@kth.se}}
\runauthor{R. Douc, E. Moulines and J. Olsson}
\affiliation{Institut T\'el\'ecom/T\'el\'ecom SudParis, Institut T\'el\'
ecom/T\'el\'ecom ParisTech\\ and KTH Royal Institute of Technology}
\address[A]{R. Douc\\
SAMOVAR, CNRS UMR 5157\\
Institut T\'el\'ecom/T\'el\'ecom SudParis\\
9 rue Charles Fourier\\
91000 Evry\\
France\\
\printead{e1}} 
\address[B]{E. Moulines\\
LTCI, CNRS UMR 8151\\
Institut T\'el\'ecom/T\'el\'ecom ParisTech\\
46, rue Barrault\\
75634 Paris Cedex 13\\
France\\
\printead{e2}}
\address[C]{J. Olsson\\
Department of Mathematics\\
KTH Royal Institute of Technology\\
SE-100 44, Stockholm\\
Sweden\\
\printead{e3}}
\end{aug}
\thankstext{both}{Supported by the Agence Nationale de la Recherche
through the 2009--2012 project Big MC.}
\thankstext{jo}{Supported by the Swedish Research Council, Grant 2011-5577.}

\received{\smonth{3} \syear{2012}}
\revised{\smonth{6} \syear{2013}}

%
\begin{abstract}
This paper discusses particle filtering in general hidden Markov models
(HMMs) and presents novel theoretical results on the long-term
stability of bootstrap-type particle filters. More specifically, we
establish that the asymptotic variance of the Monte Carlo estimates
produced by the bootstrap filter is uniformly bounded in time. On the
contrary to most previous results of this type, which in general
presuppose that the state space of the hidden state process is compact
(an assumption that is rarely satisfied in practice), our very mild
assumptions are satisfied for a large class of HMMs with possibly
noncompact state space. In addition, we derive a similar time uniform
bound on the asymptotic $\mathsf{L}^{p}$ error. Importantly, our
results hold for misspecified models; that is, we do not at all assume
that the data entering into the particle filter originate from the
model governing the dynamics of the particles or not even from an HMM.
\end{abstract}

%

\end{frontmatter}

\section{Introduction}\label{sec1}
This paper deals with estimation in general \emph
hidden Mar\-kov models (HMMs) via \emph{sequential Monte Carlo} (SMC)
\emph{methods} (or \emph{particle filters}). More specifically, we
present novel results on the numerical stability of the \emph{bootstrap
particle filter} that hold under very general and easily verifiable
assumptions. Before stating the results we provide some background.

Consider an HMM $(X_n, Y_n)_{n \in\nset}$, where the Markov chain (or
\emph{state sequence}) $(X_n)_{n \in\nset}$, taking values in some
general state space $(\set{X},
\alg{X})$, is only partially observed through the sequence $(Y_n)_{n
\in\nset}$ of \emph{observations} taking values in another general
state space $(\set{Y},
\alg{Y})$. More specifically, conditionally on the state sequence
$(X_n)_{n \in\nset}$, the observations are assumed to be conditionally
independent and such that the conditional distribution of each $Y_n$
depends on the corresponding state $X_n$ only; see, for example, \cite
{cappemoulinesryden2005} and the references therein. We denote by $\Q$
and $\chi$ the kernel and initial distribution of $(X_n)_{n \in\nset
}$, respectively. Even though $n$ is not necessarily a temporal index,
it will in the following be referred to as ``time.''

Any kind of statistical estimation in HMMs typically involves
computation of the conditional distribution of one or several hidden
states given a set of observations. Of particular interest are the \emph
{filter distributions}, where the filter distribution $\filt[\chi
]{\chunk{Y}{0}{n}}$ at time $n$ is defined as the conditional
distribution of $X_n$ given the corresponding observation history
$\chunk{Y}{0}{n} = (Y_0, \ldots, Y_n)$ (this will be our generic
notation for vectors). The problem of computing, recursively in $n$ and
in a single sweep of the data, the sequence of filter distributions is
referred to as \emph{optimal filtering}. Of similar interest are the
\emph{predictor distributions}, where the predictor distribution $\pred
[\chi]{\chunk{Y}{0}{n - 1}}$ at time $n$ is defined as the conditional
distribution of the state $X_n$ given the preceding observation history
$\chunk{Y}{0}{n - 1}$ (more precise definitions of filter and predictor
distributions are given in Section~\ref{secpreliminaries}). In this
paper we focus on the computation of these distributions, which can be
carried through in a recursive fashion according to
%
\begin{eqnarray}
\filtantri[\chi]{\chunk{Y} {0} {n}}(\set{A}) &=& \frac{\int\1_{\set{A}}(x)
\llhd{Y_n}[x] \pred[\chi]{\chunk{Y}{0}{n - 1}}(\mathrm{d} x)}{\int\llhd
{Y_n}[x]   \pred[\chi]{\chunk{Y}{0}{n - 1}}(\mathrm{d} x)},  \qquad\set{A}
\in
\alg{X},  \label{eqcorrection}
\\
\predantri[\chi]{\chunk{Y} {0} {n}}(\set{A}) &=& \int\Q(x, \set{A}) \filtantri[\chi]{
\chunk{Y} {0} {n}}(\mathrm{d} x),  \qquad\set{A} \in
\alg{X},  \label{eqprediction}
\end{eqnarray}
where $\llhd{Y_n}$ is the local likelihood of the hidden state~$X_n$
given the observation~$Y_n$. Steps \eqref{eqcorrection} and \eqref
{eqprediction} are typically referred to as \emph{correction} and \emph
{prediction}, respectively. In the correction step, the predictor $\pred
[\chi]{\chunk{Y}{0}{n - 1}}$ is, as the new observation $Y_n$ becomes
available, weighted with the local likelihood, providing the filter
$\filt[\chi]{\chunk{Y}{0}{n}}$; in the prediction step, the filter
$\filt[\chi]{\chunk{Y}{0}{n}}$ is propagated through the dynamics $\Q$
of the latent Markov chain, yielding the predictor $\pred[\chi]{\chunk
{Y}{0}{n}}$ at the consecutive time step. The correction and prediction
steps form jointly a measure-valued mapping~$\Phi$ generating
recursively the predictor distribution flow according to
\[
\predantri[\chi]{\chunk{y} {0} {n}} = \feynantri[\phi_\chi\bigl\langle\chunk{Y}{0} {n - 1} \bigr\rangle]{Y_n}
\]
(we refer again to Section~\ref{secpreliminaries} for precise definitions).

Unless the HMM is either a linear Gaussian model or a model comprising
only a finite number of possible states, exact numeric computation of
the predictor distributions is in general infeasible. Thus, one is
generally confined to using finite-dimensional approximations of these
measures, and in this paper we concentrate on the use of particle
filters for this purpose. A particle filter approximates the predictor
distribution at time $n$ by the empirical measure $\partpred[\chi
]{\chunk{Y}{0}{n - 1}}$ associated with a finite sample $(\epart
{n}{i})_{i = 1}^N$ of \emph{particles} evolving randomly and
recursively in time. At each iteration of the algorithm, the particle
sample is updated through a \emph{selection step} and a \emph{mutation
step}, corresponding directly to correction and prediction,
respectively. The selection operation duplicates or eliminates
particles with high or low posterior probability, respectively, while
the mutation operation disseminates randomly the particles in the state
space for further selection at the next iteration. The most basic
algorithm---proposed in \cite{gordonsalmondsmith1993} and referred to
as the bootstrap\vadjust{\goodbreak} particle filter---performs selection by resampling
multinomially the predictive particles $(\epart{n}{i})_{i = 1}^N$ with
probabilities proportional to the local likelihood $\wgt{n}{i} = \llhd
{Y_n}[\epart{n}{i}]$ of the particle locations; after this, the
selected particle swarm is mutated according to the dynamics $\Q$ of
the latent Markov chain. Here the self-normalized importance sampling
estimator $\partfilt{\chunk{Y}{0}{n}}$ associated with the weighted
particle sample $(\epart{n}{i}, \wgt{n}{i})_{i = 1}^N$ provides an
approximation of the filter $\filt{\chunk{Y}{0}{n}}$. Thus, subjecting
the particle sample $(\epart{n}{i})_{i = 1}^N$ to selection and
mutation is in the case of the bootstrap particle filter equivalent to
drawing, conditionally independently given $(\epart{n}{i})_{i = 1}^N$,
new particles $(\epart{n + 1}{i})_{i = 1}^N$ from the distribution
$\feyn[\phi_\chi^N \langle\chunk{Y}{0}{n - 1} \rangle]{Y_n}$ obtained
by plugging the empirical measure $\partpred[\chi]{\chunk{Y}{0}{n -
1}}$ into the filter recursion, which we denote
%
\begin{equation}
\label{eqparticlesampling}
\bigl(\epart{n + 1} {i}\bigr)_{i = 1}^N
\stackrel{\mathrm{i.i.d.}}{\sim} \feynantri[\phi_\chi^N \bigl\langle\chunk{Y}
{0} {n - 1} \bigr\rangle]{Y_n}^{\varotimes N}.
\end{equation}

Since the seminal paper \cite{gordonsalmondsmith1993}, particle filters
have been successfully applied to nonlinear filtering problems in many
different fields; we refer to the collection \cite
{doucetdefreitasgordon2001} for an introduction to particle filtering
in general and for miscellaneous examples of real-life applications.

The theory of particle filtering is an active field and there is a
number of available convergence results concerning, for example, $\Lp
{p}$ error bounds and weak convergence; see the monographs \cite
{delmoral2004,baincrisan2009} and the references therein. Most of these
results establish the convergence, as the number of particles $N$ tends
to infinity, of the particle filter for a \emph{fixed} time step $n \in
\nset$. For \emph{infinite time horizons}, that is, when $n$ tends to
infinity, convergence is less obvious. Indeed, each recursive update
\eqref{eqparticlesampling} of the particles $(\epart{n}{i})_{i=1}^N$
is based on the implicit assumption that the empirical measure
$\partpred[\chi]{\chunk{Y}{0}{n - 1}}$ associated with the ancestor
sample approximates perfectly well the predictor $\pred[\chi]{\chunk
{Y}{0}{n - 1}}$ at the previous time step; however, since the ancestor
sample is marred by an error itself, one may expect that the errors
induced at the different updating steps accumulate and, consequently,
that the total error propagated through the algorithm increases with
$n$. This would make the algorithm useless in practice. Fortunately, it
has been observed empirically by several authors (see, e.g., \cite{vanhandel2009}, Section~1.1) that the convergence of particle filters
appears to be \emph{uniform} in time also for very general HMMs.
Nevertheless, even though long-term stability is essential for the
applicability of particle filters, most existing time uniform
convergence results are obtained under assumptions that are generally
not met in practice. The aim of the present paper is thus to establish
the infinite time-horizon stability under mild and easy verifiable
assumptions, satisfied by most models for which the particle filter has
been found to be useful.

\subsection{Previous work}\label{secpreviouswork}

The first time uniform convergence result for boot\-strap-type particle
filters was obtained by Del Moral and Guionnet \cite
{delmoralguionnet2001} (see also~\cite{delmoral2004} for
refinements) using contraction properties of the mapping $\Phi$. We
recall in some detail this technique. By writing
\begin{eqnarray*}
\partpredantri[\chi]{\chunk{Y} {0} {n}} - \predantri[\chi]{\chunk{Y} {0} {n}} &=&
\underbrace{\partpredantri[\chi]{\chunk{Y} {0} {n}} - \feynantri[\phi_\chi^N
\bigl\langle\chunk{Y} {0} {n - 1} \bigr\rangle]{Y_n}}_{\mathrm{sampling\ error}}
\\
&&{} + \underbrace{\feynantri[\phi_\chi^N \bigl\langle\chunk{Y} {0}
{n - 1} \bigr\rangle ]{Y_n} - \feynantri[\phi_\chi\bigl\langle
\chunk{Y} {0} {n - 1} \bigr\rangle ]{Y_n}}_{\mathrm{initialization\ error}}
\end{eqnarray*}
one may decompose the error $\partpred[\chi]{\chunk{Y}{0}{n}} - \pred
[\chi]{\chunk{Y}{0}{n}}$ into a first error (the sampling error)
introduced by replacing $\feyn[\phi_\chi^N \langle\chunk{Y}{0}{n - 1}
\rangle]{Y_n}$ by its empirical estimate $\partpred[\chi]{\chunk
{Y}{0}{n}}$ and a second error (the initialization error) originating
from the discrepancy between empirical measure $\partpred[\chi]{\chunk
{Y}{0}{n - 1}}$ associated with the ancestor particles and the true
predictor $\pred[\chi]{\chunk{Y}{0}{n - 1}}$. The sampling error is
easy to control. One may, for example, use the Marcinkiewicz--Zygmund
inequality to bound the $\Lp{p}$ error by $c N^{-1 / 2}$, where $c \in
\rsetpos$ is a universal constant. Exponential deviation inequalities
may also be obtained. For the initialization error, we may expect that
the mapping $\feyn{Y_n}$ is in some sense contracting and thus
downscales the discrepancy between $\partpred[\chi]{\chunk{Y}{0}{n -
1}}$ and $\pred[\chi]{\chunk{Y}{0}{n - 1}}$. This is the point where
the exponential forgetting of the predictor distribution becomes
crucial. Assume, for instance, that there exists a constant $\rho\in\,
\ooint{0, 1}$ such that $\| \feyn[\mu]{\chunk{Y}{m}{n}} - \feyn[\nu
]{\chunk{Y}{m}{n}} \| \leq\rho^{n -m + 1} \| \mu- \nu\|$ for any
integers $0 \leq m \leq n$ and any probability measures $\mu$ and $\nu
$, where $\| \cdot\|$ is some suitable norm on the space of
probability measures and $\feyn{\chunk{Y}{m}{n}} \eqdef\feyn{Y_n}
\circ\feyn{Y_{n - 1}} \circ\cdots\circ\feyn{Y_m}$. Since $\feyn[\mu
]{\chunk{Y}{m}{n}}$ is the predictor distribution $\pred[\mu]{\chunk
{Y}{m}{n}}$ obtained when the hidden chain is initialized with the
distribution $\mu$ at time $m$, this means that the predictor
distribution forgets the initial distribution geometrically fast. In
addition, the forgetting rate $\rho$ is uniform with respect to the observations.
The uniformity with respect to the observations is of course the main reason why
the assumptions on the model are so stringent.

Now, decomposing similarly also the initialization error and proceeding
recursively yields the telescoping sum
%
\begin{eqnarray}\label{eqdelmoraldecomposition}
&& \partpredantri[\chi]{\chunk{Y} {0} {n}} - \predantri[\chi]{\chunk{Y} {0} {n}} \nonumber
\\
&&\qquad = \partpredantri[\chi]{\chunk{Y} {0} {n}} -
\feynantri[\phi_\chi^N \bigl\langle\chunk {Y} {0} {n - 1} \bigr\rangle]{Y_n}
\nonumber\\[-8pt]\\[-8pt]
&&\quad\qquad{}+ \sum_{k = 1}^{n - 1} \bigl( \feynantri[
\phi_\chi^N \bigl\langle\chunk {Y} {0} {k} \bigr\rangle]{
\chunk{Y} {k + 1} {n}} - \feynantri{\chunk{Y} {k + 1} {n}} \circ\feynantri[
\phi_\chi^N \bigl\langle\chunk{Y} {0} {k - 1} \bigr\rangle]{Y_k} \bigr)\nonumber
\\
&&\quad\qquad{} + \feynantri[\phi_\chi^N \bigl\langle Y_0 \bigr\rangle]{\chunk{Y} {1} {n}} - \feynantri[\phi _\chi\bigl\langle
Y_0 \bigr\rangle]{\chunk{Y} {1} {n}}.\nonumber
\end{eqnarray}
Now each term of the sum above can be viewed as a downscaling (by a
factor $\rho^{n - k}$) of the sampling error between $\partpred[\chi
]{\chunk{Y}{0}{k}}$ and $\feyn[\phi_\chi^N \langle\chunk{Y}{0}{k - 1}
\rangle]{Y_k}$ through the contraction of $\feyn{\chunk{Y}{k + 1}{n}}$.
Denoting by $\delta_n$ the $\Lp{p}$ error of $\partpred[\chi]{\chunk
{Y}{0}{n}}$ and assuming that the initial sample is obtained through
standard importance sampling, implying that $\delta_0 \leq c N^{-1 /
2}$, provides, using the contraction of $\feyn{\chunk{Y}{k + 1}{n}}$,
the uniform $\Lp{p}$ error bound $\delta_n \leq c N^{-1 / 2} \sum_{k =
0}^n \rho^{n - k} \leq c N^{-1 / 2} (1 - \rho)^{-1}$.

Even though this result is often used as a general guideline on
particle filter stability, it relies nevertheless heavily on the
assumption that the kernel $\Q$ of hidden Markov chain satisfies the
following \emph{strong mixing condition}, which is even more stringent
than the already very strong \emph{one-step global Doeblin condition}:
there exist constants $\epsilon^+ > \epsilon^- > 0$ and a probability
measure $\nu$ on $(\set{X},
\alg{X})$ such that for all $x \in\set{X}$ and $\set{A} \in
\alg{X}$,
%
\begin{equation}
\label{eqstrongmixing} \epsilon^- \nu(\set{A}) \leq\Q(x, \set{A}) \leq\epsilon^+ \nu(\set
{A}).
\end{equation}
This assumption, which in particular implies that the Markov chain is
uniformly geometrically ergodic, restricts the applicability of the
stability result in question to models where the state space $\set{X}$
is small. For Markov chains on separable metric spaces, provided that
the kernel is strongly Feller, condition \eqref{eqstrongmixing}
typically requires the state space to be compact. Some refinements have
been obtained in, for example, \cite
{leglandoudjane2003,leglandoudjane2004,delmoral2004,oudjanerubenthaler2005,tadicdoucet2005,cappemoulinesryden2005,olssonryden2008,doucmoulinesolsson2008,crisanheine2008,heinecrisan2008}.

The long-term stability of particle filters is also related to the
boundedness of the asymptotic variance. The first central limit theorem
(CLT) for bootstrap-type particle filters was derived by Del Moral and
Guionnet \cite{delmoralguionnet1999}. More\vspace*{1.5pt} specifically, it was shown
that the normalized Monte Carlo error $\sqrt{N}(\partpred[\chi]{\chunk
{Y}{0}{n - 1}} h - \pred[\chi]{\chunk{Y}{0}{n - 1}} h)$ tends weakly,
for a fixed $n \in\nsetpos$ and as the particle population size $N$
tends to infinity, to a zero mean normal-distributed variable with
variance $\variance[2]{\chunk{Y}{0}{n - 1}}{h}$. Here we have used the
notation $\mu h \eqdef\int h(x)   \mu(\mathrm{d} x)$ to denote
expectations. The original proof of the CLT was later simplified and
extended to more general particle filtering algorithms in \cite
{kuensch2005,chopin2002,doucmoulines2008,doucmoulinesolsson2008,doucetjohansen2008};
in Section~\ref{secpreliminaries} we recall in detail the version
obtained in \cite{doucmoulines2008} and provide an explicit expression
of the asymptotic variance $\variance[2]{\chunk{Y}{0}{n - 1}}{h}$. As
shown first by \cite{delmoralguionnet2001}, Theorem~3.1, it is
possible, using the strong mixing assumption described above, to bound
uniformly also the asymptotic variance $\variance[2]{\chunk{Y}{0}{n -
1}}{h}$ by similar forgetting-based arguments. Here a key ingredient is
that the particles $(\epart{n}{i})_{i=1}^N$ obtained at the different
time steps become, asymptotically as $N$ tends to infinity,
statistically\vspace*{1pt} independent. Consequently, the total asymptotic variance
of $\sqrt{N}(\partpred[\chi]{\chunk{Y}{0}{n - 1}} h - \pred[\chi]{\chunk
{Y}{0}{n - 1}} h)$ is obtained by simply summing up the asymptotic
variances of the error terms $\sqrt{N}(\feyn[ \phi_\chi^N \langle\chunk
{Y}{0}{k} \rangle]{\chunk{Y}{k + 1}{n}}h - \feyn{\chunk{Y}{k + 1}{n}}
\circ\feyn[ \phi_\chi^N \langle\chunk{Y}{0}{k - 1} \rangle]{Y_k}h)$
in the decomposition~\eqref{eqdelmoraldecomposition}. Finally,
applying again the contraction of the composed mapping $\feyn{\chunk
{Y}{m}{n}}$ yields a uniform bound on the total asymptotic variance in
accordance with the calculation above. In \cite
{doucgariviermoulinesolsson2009}, a similar stability result was
obtained for a particle-based version of the \emph{forward-filtering
backward-simulation algorithm} (proposed in \cite
{godsilldoucetwest2004}); nevertheless, also the analysis of this work
relies completely on the assumption of strong mixing of the latent
Markov chain, which, as we already pointed out, does not hold for most
models used in practice.

A first breakthrough toward stability results for noncompact state
spaces was made in \cite{vanhandel2009}. This work establishes, again
for bootstrap-type particle filters, a~uniform time average convergence
result of form
%
\begin{equation}
\label{eqvanhandelconvergence}
\lim_{N \to\infty} \sup_{n \in\nset}\espBigg{n^{-1} \sum_{k = 1}^n \bigl\|
\partfiltantri[\chi]{\chunk{Y} {0} {k}} - \filtantri[\chi]{\chunk{Y} {0} {k}} \bigr\|
_{\mathrm{BL}}} = 0,
\end{equation}
where $\| \cdot\|_{\mathrm{BL}}$ denotes the dual bounded-Lipschitz
norm. This result, obtained as a special case of a general
approximation theorem derived in the same paper, was established under
very weak assumptions on the local likelihood (supposed to be bounded
and continuous) and the Markov kernel (supposed to be Feller). These
assumptions are, together with the basic assumption that the hidden
Markov chain is positive Harris and aperiodic, satisfied for a large
class of HMMs with possibly noncompact state spaces. Nevertheless, the
proof is heavily based on the assumption that the particles evolve
according to exactly the same model dynamics as the observations
entered into the algorithm, in other words, that the model is perfectly
specified. This is of course never true in practice. In addition, the
convergence result \eqref{eqvanhandelconvergence} does not, in
contrast to $\Lp{p}$ bounds and CLTs, provide a rate of convergence of
the algorithm.



\subsection{Approach of this paper}\label{secthispaper}

In this paper we return to more standard convergence modes and
reconsider the asymptotic variance and $\Lp{p}$ error of bootstrap
particle filters. As noticed by Johansen and Doucet \cite
{doucetjohansen2008}, restricting the analysis of bootstrap-type
particle filters does not imply a significant loss of generality, as
the CLT for more general \emph{auxiliary particle filters} \cite
{pittshephard1999} can be straightforwardly obtained by applying the
bootstrap filter CLT to a somewhat modified HMM incorporating the
so-called \emph{adjustment multiplier weights} of the auxiliary
particle filter into the model dynamics. Our aim is to establish that
the asymptotic variance and $\Lp{p}$ error are stochastically bounded
in the noncompact case. Recall that a sequence $(\mu_n)_{n \in\nset}$
of probability measures on $(\rset, \mathcal{B}(\rset))$ is \emph
{tight} if for all $\epsilon> 0$ there exists a compact interval $\set
{I} = [- a, a] \subset\rset$ such that $\mu_n(\set{I}^c) \leq\epsilon
$ for all $n$. In addition, we call a sequence $(Z_n)_{n \in\nset}$ of
random variables, with $Z_n \sim\mu_n$, tight if the sequence $(\mu
_n)_{n \in\nset}$ of marginal\vspace*{1pt} distributions is tight. In this paper,
we show that the sequence $(\sigma_\chi^2 \langle\chunk{Y}{0}{n - 1}
\rangle(h))_{n \in\nsetpos}$ of asymptotic variances is tight for
\emph{any stationary and ergodic sequence} $(Y_n)_{n \in\nset}$ \emph
{of observations}. In particular, we do not at all assume that the
observations originate from the model governing the dynamics of the
particle filter or not even from an HMM.

Our proofs are based on novel coupling techniques developed in \cite
{doucmoulines2011} (and going back to \cite{kleptsynaveretennikov2008}
and \cite{doucfortmoulinespriouret2009}) with the purpose of
establishing the convergence of the \emph{relative entropy} for
misspecified HMMs. In our analysis, the strong mixing assumption~\eqref
{eqstrongmixing} is replaced by the considerably weaker $r$-\emph{local Doeblin condition}~\eqref{eqdefinition-LD-set}. This
assumption is, for instance, trivially satisfied (for $r = 1$) if there
exist a measurable set $\set{C} \subseteq\set{X}$, a probability
measure $\lambda_{\set{C}}$ on $(\set{X},
\alg{X})$ such that $\lambda_{\set{C}}(\set{C}) = 1$ and positive
constants $0 < \epsilon^-_{\set{C}} < \epsilon^+_{\set{C}}$ such that
for all $x \in\set{X}$ and all $\set{A} \in
\alg{X}$,
%
\begin{equation}
\label{eq1-localDoeblin} \epsilon^-_{\set{C}} \lambda_{\set{C}}(\set{A}) \leq
\Q(x, \set{A} \cap \set{C}) \leq\epsilon^+_{\set{C}} \lambda_{\set{C}}(
\set{A}),
\end{equation}
a condition that is easily verified for many HMMs with noncompact
state space [we emphasize, however, that assumption \eqref
{eqdefinition-LD-set} is even weaker than~\eqref{eq1-localDoeblin}].

%
\begin{remark}
Finally, we remark that Del Moral and Guionnet \cite
{delmoralguionnet2001} studied the stability of SMC methods within the
framework of a general \emph{normalized Feynman--Kac prediction model}
consisting of a sequence $(\fm{n})_{n \in\nset}$ of measures-defined
recursively on a sequence $(\set{E}_n, \field{E}_n)_{n \in\nset}$ of
measurable spaces by
\[
\fm{n + 1}(\set{A}) = \frac{\int\pot[x]{n} \K_n(x, \set{A})   \fm
{n}(\mathrm{d} x)}{\int\pot[x]{n}  \fm{n}(\mathrm{d} x)},  \qquad\set{A} \in \field{E}_{n + 1}, %
\]
where $\pot{n}$ is a positive potential function on $\set{E}_n$, and $\K
_n$ is a Markov transition kernel from $(\set{E}_n, \field{E}_n)$ to
$(\set{E}_{n + 1}, \field{E}_{n + 1})$; see also \cite{delmoral2004}, Section~2.3. Conditionally on the observations
$(Y_n)_{n \in\nset}$, the flow of predictor distributions can
obviously, by \eqref{eqcorrection}--\eqref{eqprediction}, be
formulated as a normalized Feynman--Kac prediction model by letting
$(\set{E}_n, \field{E}_n) \equiv(\set{X}, \field{X})$, $\fm{n} \equiv
\pred[\chi]{\chunk{Y}{0}{n - 1}}$, $\pot{n} \equiv\llhd{Y_n}$, and $\K
_n \equiv\Q$, $n \in\nset$. Imposing, as in \cite
{delmoralguionnet2001}, the assumption that the transition kernels $(\K
_n)_{n \in\nset}$ satisfy jointly the global Doeblin condition \eqref
{eqstrongmixing} provides a mixing rate $\rho$ that is uniform in the
observations, and any stability result obtained for fixed observations
holds thus automatically true also when the observations are allowed to
vary randomly.

Similarly, we could in the present paper have taken directly the
recursion \eqref{eqcorrection}--\eqref{eqprediction} as a starting
point, by suppressing its connection with HMMs and describing the same
as a normalized Feynman--Kac prediction model indexed by a stationary
and ergodic sequence $(Y_n)_{n \in\nset}$ of random parameters.
However, as the results obtained in \cite{doucmoulines2011}, which are
fundamental in our analysis, describes the convergence of the
normalized log-likelihood function for general HMMs, a quantity whose
interpretation is not equally clear in the context of Feynman--Kac
models, we have chosen to express our results in the language of HMMs
as well.
\end{remark}

To sum up, the contribution of the present paper is twofold, since:
\begin{itemize}
\item we present time uniform bounds providing also the rate of
convergence in $N$ of the particle filter for HMMs with possibly
noncompact state space;
\item we establish long-term stability of the particle filter also in
the case of misspecification, that is, when the stationary law of the
observations entering the particle filter differs from that of the HMM
governing the dynamics of the particles~$(\epart{n}{i})_{i=1}^N$.
\end{itemize}

\subsection{Outline of the paper}
The paper is organized as follows. Section~\ref{secpreliminaries}
provides the main notation and definitions. It also introduces the
concepts of HMMs and bootstrap particle filters. In Section~\ref
{secmainresults} our main results are stated together with the main
layouts of the proofs. Section~\ref{secapplications} treats some
examples, and Section~\ref{secproofs} and Section~\ref{sectechlemmas}
provide the full details of our proofs.

\section{Preliminaries}\label{secpreliminaries}
\subsection{Notation}

We preface the introduction of HMMs with some notation. Let $(\set{X},
\field{X})$ be a measurable space, where $\field{X}$ is a countably
generated $\sigma$-field. Denote by $\bmf{
\alg{X}}$ [resp., $\pmf{
\alg{X}}$] the space of bounded (resp., bounded and nonnegative) $\field
{X}/\borel{\rset}$-measurable functions on $\set{X}$ equipped with the
supremum norm $\| f \|_\infty\eqdef\sup_{x \in\set{X}} |f(x)|$. In
addition, denote by $\probmeas{
\alg{X}}$ the set of probability measures on $(\set{X},
\alg{X})$. Let $\K\dvtx  \set{X} \times
\alg{X} \rightarrow\rset_+$ be a finite kernel on $(\set{X}, \field
{X})$, that is, for each $x \in\set{X}$, the mapping $\K(x, \cdot) \dvtx
\alg{X} \ni\set{A} \mapsto\K(x, \set{A})$ is a finite measure on
$\field{X}$, and for each $\set{A} \in\field{X}$, the function $\K(\cdot, \set{A}) \dvtx  \set{X} \ni x \mapsto\K(x, \set{A})$ is $\field{X}/\borel
{\ccint{0,1}}$-measurable. If $\K(x, \cdot)$ is a probability measure
on $(\set{X}, \field{X})$ for all $x \in\set{X}$, then the kernel $\K$
is said to be Markov. A kernel induces two integral operators, the
first acting on the space $\meas{\field{X}}$ of $\sigma$-finite
measures on $(\set{X},
\alg{X})$ and the other on $\bmf{\field{X}}$. More specifically, for
$\mu\in\meas{\field{X}}$ and $f \in\bmf{\field{X}}$, we define the measure
\[
\mu\K \dvtx
\alg{X} \ni\set{A} \mapsto\int\K(x, \set{A}) \mu(\mathrm{d} x)
\]
and the function
\[
\K f \dvtx  \set{X} \ni x \mapsto\int f\bigl(x'\bigr) \K\bigl(x, \mathrm{d}
x'\bigr).  %
\]
Moreover, the \emph{composition} (or \emph{product}) of two kernels $\K
$ and $\M$ on $(\set{X}, \field{X})$ is defined as
\[
\K\M\dvtx  \set{X} \times
\alg{X} \ni(x, \set{A}) \mapsto \int\M
\bigl(x', \set{A}\bigr) \K\bigl(x,\mathrm{d} x'\bigr).
\]

\subsection{Hidden Markov models}
Let $(\set{X},
\alg{X})$ and $(\set{Y},
\alg{Y})$ be two measurable spaces. We specify the HMM as follows. Let
$\Q\dvtx  \set{X} \times
\alg{X} \rightarrow[0, 1]$ and $\G\dvtx  \set{X} \times
\alg{Y} \rightarrow[0, 1]$ be given Markov kernels, and let $\chi$ be
a given initial distribution on $(\set{X},
\alg{X})$. In this setting, define the Markov kernel
%
\begin{eqnarray}
\T\bigl((x, y), \set{A}\bigr) \eqdef\iint\1_{\set{A}}\bigl(x',
y'\bigr) \Q\bigl(x, \mathrm{d} x'\bigr) \G
\bigl(x', \mathrm{d} y'\bigr),\nonumber
\\
\eqntext{(x, y) \in\set{X} \times\set{Y},  \set{A} \in
\alg{X}
\algprod
\alg{Y},}
\end{eqnarray}
on the product space $(\set{X} \times\set{Y},
\alg{X} \algprod
\alg{Y})$. Let $(X_n, Y_n)_{n \in\nset}$ be the canonical Markov chain
induced by $\T$ and the initial distribution $
\alg{X} \algprod
\alg{Y} \ni\set{A} \mapsto\break \int\1_{\set{A}}(x, y)   \chi(\mathrm{d} x)
\G(x, \mathrm{d} y)$. The bivariate process $(X_n, Y_n)_{n \in\nset}$ is
what we refer to as the HMM. We shall denote by $\probcan$ and $\espcan
$ the probability measure and corresponding expectation associated with
the HMM on the canonical space $((\set{X} \times\set{Y})^\nset, (
\alg{X} \algprod
\alg{Y})^{\algprod\nset})$. We assume that the observation kernel $\G$
is nondegenerated in the sense that there exist a $\sigma$-finite
measure $\yref$ on $(\set{Y},
\alg{Y})$ and a measurable function $g \dvtx  \set{X} \times\set{Y}
\rightarrow\,\ooint{0, \infty}$ such that
\[
\G(x, \set{A}) = \int\1_{\set{A}}(y) g(x, y) \yref(\mathrm{d} y),  \qquad x \in
\set{X},  \set{A} \in
\alg{Y}.  %
\]
For a given observation $y \in\set{Y}$, we let
\[
\llhd{y} \dvtx  \set{X} \ni x \mapsto g(x, y)
\]
denote the local likelihood function of the state given the
corresponding observation $y$.

When operating on HMMs we are in general interested in computing
expectations of type $\espcan(h( \chunk{\Xcan}{k}{\ell}) | \chunk{\Ycan
}{0}{m})$ for integers $(k, \ell, m) \in\nset^3$ with $k \leq\ell$
and functions $h \in\bmf{\set{X}^{\ell- k + 1}}$. Of particular
interest are quantities of form $\espcan(h(\Xcan_n) | \chunk{\Ycan
}{0}{n - 1})$ or $\espcan(h(\Xcan_n) | \chunk{\Ycan}{0}{n})$, and the
term optimal filtering refers to problem of computing, recursively in
$n$, such conditional distributions and expectations as new data
becomes available. For any record $\chunk{y}{k}{m} \in\set{Y}^{m - k +
1}$ of observations, let $\Lkh{\chunk{y}{k}{m}}$ be the unnormalized
kernel on $(\set{X}, \field{X})$ defined by
%
\begin{eqnarray}\label{eqdefinition-Lkh}
\Lkhantri{\chunk{y} {k} {m}}(x_k, \set{A}) \eqdef
\idotsint \1_{\set
{A}}(x_{m + 1}) \prod
_{\ell= k}^m \llhd{y_\ell}[x_\ell]
\Q(x_\ell, \mathrm{d} x_{\ell+1}),
\nonumber\\[-10pt]\\[-10pt]
\eqntext{x_k \in\set{X},  \set{A} \in
\alg{X},}
\end{eqnarray}
with the convention
%
\begin{equation}
\label{eqconvention} \Lkhantri{\chunk{y} {k} {m}}(x, \set{A}) \eqdef\delta_x(
\set{A})\qquad\mbox {for }k > m
\end{equation}
(where $\delta_x$ denotes the Dirac mass at point $x$). Note that the
function $\chunk{y}{0}{n - 1} \mapsto\chi\Lkh{\chunk{y}{0}{n-1}} \1
_\set{X}$ is exactly the density of the observations $\chunk{Y}{0}{n -
1}$ (i.e., the likelihood function) with respect to $\yref^{\algprod
n}$. Also note that for any $\ell\in\{k, \ldots, m - 1\}$,
%
\begin{equation}
\label{eqLkh-composition} \Lkhantri{\chunk{y} {k} {m}} = \Lkhantri{\chunk{y} {k} {\ell}} \Lkhantri{\chunk{y} {\ell+ 1} {m}}.
\end{equation}
Let $\pred[\chi]{\chunk{y}{k}{m}}$ be the probability measure defined by
%
\begin{equation}
\label{eqexpressionpred} \predantri[\chi]{\chunk{y} {k} {m}} (\set{A}) \eqdef\frac{\chi\Lkh{\chunk
{y}{k}{m}} \1_\set{A} }{\chi\Lkh{\chunk{y}{k}{m}} \1_\set{X}}, \qquad \set{A} \in
\alg{X}.
\end{equation}
Note that this implies that $\pred[\chi]{\chunk{y}{k}{m}} = \chi$ when
$k > m$. Using the notation, it can be shown (see, e.g., \cite{cappemoulinesryden2005}, Proposition~3.1.4) that for any $h \in\bmf
{\field{X}}$,
\[
\espcan \bigl( h(X_n) \midbar \chunk{\Ycan} {0} {n - 1} \bigr) = \int
h(x) \predantri[\chi]{\chunk{\Ycan} {0} {n - 1}}(\mathrm{d} x),  %
\]
that is, $\pred[\chi]{\chunk{\Ycan}{0}{n - 1}}$ is the predictor of
$X_n$ given the observations $\chunk{\Ycan}{0}{n - 1}$. From definition
\eqref{eqexpressionpred} one immediately obtains the recursion
\[
\label{eqfilterrecursion} \predantri[\chi]{\chunk{y} {0} {n}}(\set{A}) = \frac{\pred[\chi]{\chunk
{y}{0}{n - 1}} \Lkh{y_n} \1_\set{A}}{\pred[\chi]{\chunk{y}{0}{n - 1}}
\Lkh{y_n} \1_\set{X}} =
\frac{\int\llhd{y_n}[x] \Q(x, \set{A})   \pred[\chi]{\chunk
{y}{0}{n - 1}}(\mathrm{d} x)}{\int\llhd{y_n}[x]   \pred[\chi]{\chunk
{y}{0}{n - 1}}(\mathrm{d} x)},  \quad\set{A} \in
\alg{X},
\]
which can be expressed in condensed form as
%
\begin{equation}
\label{eqfilterrecursioncondensed} \predantri[\chi]{\chunk{y} {0} {n}} = \feyn{y_n}\bigl(
\predantri[\chi]{\chunk{y} {0} {n - 1}}\bigr)
\end{equation}
with $\feyn{y_n}$ being the measure-valued transformation
\[
\feyn{y_n}\dvtx  \probmeas{\field{X}} \ni\mu\mapsto\upstep{y_n}(
\mu) \Q
\]
and $\upstep{y_n}$ transforms a measure $\mu\in\probmeas{\field{X}}$
into the measure
\[
\upstep{y_n}(\mu) \dvtx
\alg{X} \ni\set{A} \mapsto
\frac{\int\1_{\set{A}}(x) \llhd{y_n}[x]
\mu(\mathrm{d} x)}{\int \llhd{y_n}[x]   \mu(\mathrm{d} x)} %
\]
in $\probmeas{\field{X}}$. By introducing the filter distributions
\[
\filtantri[\chi]{\chunk{y} {0} {n}} \eqdef\upstep{y_n}\bigl( \predantri[\chi]{
\chunk {y} {0} {n - 1}} \bigr),  %
\]
satisfying, for all $h \in\bmf{\field{X}}$,
\[
\espcan \bigl( h(X_n) \midbar\chunk{\Ycan} {0} {n} \bigr) = \int
h(x) \predantri[\chi]{\chunk{\Ycan} {0} {n}}(\mathrm{d} x) %
\]
(see again \cite{cappemoulinesryden2005}, Proposition~3.1.4), we may
express one iteration of the filter recursion in terms of the two operations
\[
\predantri[\chi]{\chunk{y} {0} {n - 1}} \mathop{-\!\!\!-\!\!\!-\!\!\!-\!\!\!\rightarrow}_{\mathrm{Updating}}^{\upstep
{y_n}} \filtantri[\chi]{\chunk{y} {0} {n}} \mathop{-\!\!\!-\!\!\!-\!\!\!-\!\!\!-\!\!\!\rightarrow}_{\mathrm{Prediction}}^{\Q } \predantri[\chi]{\chunk{y} {0} {n}}.  %
\]

As mentioned in the \hyperref[sec1]{Introduction}, it is in general infeasible to obtain
closed-form solutions to the recursion \eqref
{eqfilterrecursioncondensed}. In the following section we discuss
how approximate solutions to \eqref{eqfilterrecursioncondensed} can
be obtained using particle filters, with focus set on the bootstrap
particle filter proposed in \cite{gordonsalmondsmith1993}.

\subsection{The bootstrap particle filter}

In the following we assume that all random variables are defined on a
common probability space $(\Omega,
\alg{A}, \probalone)$. The bootstrap particle filter updates
sequentially a set of weighted simulations in order to approximate
online the flow of predictor and filter distributions. In order to
describe precisely how this is done for a given sequence $(y_n)_{n \in
\nset}$ of observations, we proceed inductively and assume that we are
given a sample of $\set{X}$-valued random draws $(\epart{n}{i})_{i =
1}^N$ (the particles) such that the empirical measure
\[
\label{eqSMCapprox} \partpredantri[\chi]{\chunk{y} {0} {n-1}} \eqdef\frac{1}{N} \sum
_{i = 1}^N \delta_{\epart{n}{i}}
\]
associated with these draws \emph{targets} the predictor $\pred[\chi
]{\chunk{y}{0}{n-1}}$ in the sense that $\partpred[\chi]{\chunk{y}{0}{n
- 1}} h = \sum_{i = 1}^N h(\epart{n}{i}) / N$ estimates $\pred[\chi
]{\chunk{y}{0}{n - 1}}h$ for any $h \in\bmf{\field{X}}$. In order to
form a new particle sample $(\epart{n + 1}{i})_{i = 1}^N$ approximating
the predictor $\pred[\chi]{\chunk{y}{0}{n}}$ at the subsequent time
step, we replace, in \eqref{eqfilterrecursioncondensed}, the true
predictor $\pred[\chi]{\chunk{y}{0}{n - 1}}$ by the particle estimate
$\partpred[\chi]{\chunk{y}{0}{n - 1}}$, and pass the latter through the
updating and prediction steps. This yields, after updating, the
self-normalized approximation
%
\begin{equation}
\label{eqparticleapproximationupdatingstep} \partfiltantri[\chi]{\chunk{y} {0} {n}} \eqdef\upstep{y_n}
\bigl(\partpredantri[\chi ]{\chunk{y} {0} {n - 1}}\bigr) = \sum
_{i = 1}^N \frac{\wgt{n}{i}}{\wgtsum{n}} \delta_{\epart{n}{i}}
\end{equation}
of the filter $\filt[\chi]{\chunk{y}{0}{n}}$,
where we have introduced the \emph{importance weights} $\wgt{n}{i}
\eqdef\llhd{y_n}[\epart{n}{i}]$, $i \in\{1, 2, \ldots, N \}$ and
$\wgtsum{n} \eqdef\sum_{i = 1}^N \wgt{n}{i}$. Moreover, by propagating
filter approximation \eqref{eqparticleapproximationupdatingstep}
through prediction step one obtains the approximation
%
\begin{equation}
\label{eqpartmixture} \sum_{i = 1}^N
\frac{\wgt{n}{i}}{\wgtsum{n}} \Q\bigl(\epart{n} {i}, \cdot\bigr)
\end{equation}
of the predictor $\pred[\chi]{\chunk{y}{0}{n}}$. Finally, the sample
$(\epart{n + 1}{i})_{i = 1}^N$ is generated by simulating $N$
conditionally independent draws from the mixture in \eqref
{eqpartmixture} using the following algorithm:\vspace*{6pt}
\begin{algorithmic}
\STATE set $\wgtsum{n} \gets0$
\FOR{$i = 1 \to N$}
\STATE set $\wgt{n}{i} \gets\llhd{y_n}[\epart{n}{i}]$\vspace*{2pt}
\STATE set $\wgtsum{n} \gets\wgtsum{n} + \wgt{n}{i}$
\ENDFOR
\FOR{$i = 1 \to N$}
\STATE draw $\ind{n}{i} \sim(\wgt{n}{\ell} / \wgtsum{n})_{\ell= 1}^N$
\STATE draw $\epart{n + 1}{i} \sim\Q(\epart{n}{\ind{n}{i}}, \cdot)$
\ENDFOR
\end{algorithmic}\vspace*{6pt}
After this, the empirical measure $\partpred[\chi]{\chunk{y}{0}{n}} =
\sum_{i = 1}^N \delta_{\epart{n + 1}{i}} / N$ is returned as an
approximation of $\pred[\chi]{\chunk{y}{0}{n}}$. In the scheme above,
the operation $\sim$ means implicitly that all draws (for different
$i$'s) are conditionally independent. Moreover, the operation $\ind
{n}{i} \sim(\wgt{n}{\ell} / \wgtsum{n})_{\ell= 1}^N$ means that each
index $\ind{n}{i}$ is simulated according to the discrete probability
distribution generated by the normalized importance weights $(\wgt
{n}{\ell} / \wgtsum{n})_{\ell= 1}^N$. The procedure described above is
repeated recursively in order to produce particle approximations of the
predictor and filter distributions at all time steps. The algorithm is
typically initialized by drawing $N$ i.i.d. particles $(\epart{0}{i})_{i
= 1}^N$ from the initial distribution $\chi$ and letting $\sum_{i =
1}^N \delta_{\epart{0}{i}} / N$ be an estimate of $\chi$.

As mentioned in the \hyperref[sec1]{Introduction}, the asymptotic properties, as the
number $N$ of particles tends to infinity, of the bootstrap particle
filter output are well investigated. When it concerns weak convergence,
Del Moral and Guionnet \cite{delmoralguionnet1999} established the
following CLT. Define for $h \in\bmf{\field{X}}$,
%
\begin{equation}
\label{eqexpression-asymptotic-variance}
 \varianceantri[2]{\chunk{y} {0} {n - 1}} {h} \eqdef\sum
_{k = 0}^n \predantri[\chi]{\chunk{y} {0} {k - 1}} \biggl(
\frac{\Lkh{\chunk{y}{k}{n - 1}} h - \pred[\chi]{\chunk{y}{0}{n
- 1}} h \times\Lkh{\chunk{y}{k}{n - 1}} \1_\set{X}} {
\pred[\chi]{\chunk{y}{0}{k - 1}} \Lkh{\chunk{y}{k}{n - 1}} \1_\set
{X}} \biggr)^2.\hspace*{-25pt}
\end{equation}

%
\begin{theorem}[{\cite{delmoralguionnet1999}}] \label{thmcltparticle}
For all $h \in\bmf{\field{X}}$ and all $\chunk{y}{0}{n - 1} \in\set
{Y}^n$ such that $\supnormTxt{\llhd{y_\ell}} < \infty$ for all $y_\ell
$, it holds, as $N \to\infty$,%
\begin{equation}
\label{eqdlim-clt} \sqrt{N} \bigl(\partpredantri[\chi]{\chunk{y} {0} {n - 1}}h - \predantri[
\chi]{\chunk {y} {0} {n - 1}}h\bigr) \dlim\varianceantri{\chunk{y} {0} {n - 1}} {h} Z,
\end{equation}
where $\variance{\chunk{y}{0}{n - 1}}{h}$ is defined in \eqref
{eqexpression-asymptotic-variance}, and $Z$ is a standard
normal-distributed random variable.
\end{theorem}

The next corollary states the corresponding CLT for the particle
filter. Also this result is standard and is an immediate consequence of
Theorem~\ref{thmcltparticle}, the fact that $\wgtsum{n} / N$
converges, for all $\chunk{y}{0}{n} \in\set{Y}^{n + 1}$ and as $N$
tends to infinity, in probability to $\pred[\chi]{\chunk{y}{0}{n - 1}}
\llhd{y_n}$ (see, e.g., \cite{doucmoulines2008}, Theorem~10), and
Slutsky's lemma. Let again $\sigma^2_\chi$ be given by \eqref
{eqexpression-asymptotic-variance}, and define for $\chunk{y}{0}{n}
\in\set{Y}^{n + 1}$ and $h \in\bmf{\field{X}}$,
%
\begin{equation}
\label{eqdeffiltvariance} \filtvarianceantri[2]{\chunk{y} {0} {n}} {h} \eqdef\frac{\variance[2]{\chunk
{y}{0}{n - 1}}{\llhd{y_n} \{ h - \filt[\chi]{\chunk{y}{0}{n}}h \}
}}{(\pred[\chi]{\chunk{y}{0}{n - 1}} \llhd{y_n})^2}.
\end{equation}

%
\begin{corollary} \label{thmcltparticlefilter}
For all $h \in\bmf{\field{X}}$ and $\chunk{y}{0}{n} \in\set{Y}^{n +
1}$ such that $\supnormTxt{\llhd{y_\ell}} < \infty$ for all $y_\ell$ it
holds, as $N \to\infty$,
%
\begin{equation}
\label{eqdlim-cltparticlefilter} \sqrt{N} \bigl(\partfiltantri[\chi]{\chunk{y} {0} {n}}h - \filtantri[\chi]{
\chunk {y} {0} {n}}h\bigr) \dlim\filtvarianceantri{\chunk{y} {0} {n}} {h} Z,
\end{equation}
where $\filtvariance{\chunk{y}{0}{n}}{h}$ is defined in \eqref
{eqdeffiltvariance} and $Z$ is a standard normal-distributed random variable.
\end{corollary}

When the observations $(Y_n)_{n \in\nset}$ entering the particle
filter are random, the sequences $(\variance[2]{\chunk{Y}{0}{n -
1}}{h})_{n \in\nset}$ and $(\filtvariance[2]{\chunk{Y}{0}{n}}{h})_{n
\in\nset}$ of asymptotic variances are $(\Yfilt{n})_{n \in\nset
}$-adapted stochastic processes, where $(\Yfilt{n})_{n \in\nset}$ is
the natural filtration of the observation process. The aim of the next
section is to establish that these sequences are tight. \emph{Importantly},
\emph{we assume in the following that the observations
$(Y_n)_{n \in\nset}$ entering the particle filter algorithm is an
arbitrary $\probalone$-stationary sequence taking values in $\set{Y}$}.
The stationary process $(Y_n)_{n \in\nset}$ can be embedded into a
stationary process $(Y_n)_{n \in\zset}$ with doubly infinite time. In
particular, we do not at all assume that the observations originate
from the model governing the dynamics of the particles; indeed, in the
framework we consider, we do not even assume that the observations
originate from an HMM.

\section{Main results and assumptions}
\label{secmainresults}Before listing our main assumptions, we recall
the definition of a $r$-local Doeblin set.
%
\begin{definition} \label{defilocal-Doeblin}
Let $r \in\nsetpos$. A set $\set{C} \in\field{X}$ is $r$-\emph{local
Doeblin} with respect to $\{\Q, g\}$ if there exist positive functions $\epsilon
^-_{\set{C}}\dvtx  \set{Y}^r \to\rset^+$ and $\epsilon^+_{\set{C}}\dvtx  \set
{Y}^r \to\rset^+$, a family $\{ \probdoeblin{\set{C}}{z}; z \in\set
{Y}^r \}$ of probability measures, and a family $\{\varphi_{\set{C}}
\langle z \rangle; z \in\set{Y}^r \}$ of positive functions such that
for all $z \in\set{Y}^r$, $\probdoeblin{\set{C}}{z}(\set{C}) =1$ and
for all $\set{A} \in\field{X}$ and $x \in\set{C}$,
%
\begin{eqnarray}
\label{eqdefinition-LD-set} \epsilon^-_{\set{C}} \langle z \rangle
\varphi_{\set{C}} \langle z \rangle(x) \probdoeblin{\set{C}} {z} (\set{A})
&\leq&\Pblock{z}(x, \set {A} \cap\set{C})
\nonumber\\[-8pt]\\[-8pt]
&\leq& \epsilon^+_{\set{C}} \langle z\rangle \varphi _{\set{C}} \langle z \rangle(x) \probdoeblin{\set{C}} {z}(
\set{A}).\nonumber
\end{eqnarray}
\end{definition}
\begin{longlist}[(A1)]\label{assum:likelihoodDrift}
\item[(A1)] The process $\sequence{Y}$ is strictly stationary and
ergodic. Moreover, there exist an integer $r \in\nsetpos$ and a set
$\set{K} \in\field{Y}^{\varotimes r}$ such that the following holds:
\begin{longlist}[(iii)]
\item[(i)] The process $\sequence{Z}$, where $Z_n \eqdef\chunk{Y}{n r}{(n + 1)r -
1}$, is ergodic and such that $
\prob{Z_0 \in\set{K}} > 2/3$.
\item[(ii)] For\vspace*{1pt} all $\eta> 0$ there exists an $r$-local Doeblin set $\set{C}
\in\field{X}$ such that for all \mbox{$\chunk{y}{0}{r-1} \in\set{K}$},
%
\begin{equation}
\label{eqbound-eta-G22} \sup_{x \in\set{C}^c} \Pblockantri{\chunk{y} {0} {r-1}}(x,
\set{X}) \leq\eta \supnormantri{\Pblockantri{\chunk{y} {0} {r-1}}(\cdot,\set{X})} < \infty
\end{equation}
and
%
\begin{equation}
\label{eqlower-bound} \inf_{\chunk{y}{0}{r - 1} \in\set{K}} \frac{\epsilon_{\set{C}}^-
\langle\chunk{y}{0}{r - 1} \rangle}{\epsilon_{\set{C}}^+ \langle\chunk
{y}{0}{r-1} \rangle} > 0,
\end{equation}
where the functions $\epsilon^+_{\set{C}}$ and $\epsilon^-_{\set{C}}$
are given in Definition~\ref{defilocal-Doeblin}.
\item[(iii)]\label{itemcondition-minoration}
There exists a set $\set{D} \in\field{X}$ such that
%
\begin{equation}
\label{eqcondition-minoration} \espBig{\ln^- \inf_{x \in\set{D}} \delta_x
\Pblockantri{\chunk{Y} {0} {r-1}} \1 _{\set{D}}} < \infty.
\end{equation}
\end{longlist}

\item[(A2)]\label{assum:majo-g}
(i)\label{itemassum-pos-g} $g(x, y) > 0$ for all $(x, y) \in
\set{X} \times\set{Y}$.

(ii)\label{itemassum-supg} $\esp{\ln^+ \supnorm{\llhd{Y_0}}} <
\infty$.
\end{longlist}

%
\begin{remark}
Assumption \hyperref[assum:likelihoodDrift]{(A1)}(i) is
inherited from \cite{doucmoulines2011}. To get some rationale behind
the constant $2/3$ appearing in the assumption, note that the same is
in fact equivalent to $1 -
\prob{Z_0 \in\set{K}} < 2
\prob{Z_0 \in\set{K}} - 1$. In that case there exist $0 < \gamma^- <
\gamma^+$ such that
\[
1 -
\prob{Z_0 \in\set{K}} < \gamma^- < \gamma^+ <
2
\prob{Z_0 \in\set{K}} - 1,  %
\]
which is equivalent to $
\prob{Z_0 \in\set{K}} > \max\{ 1-\gamma^-,(1+\gamma^+)/2 \}$. The
latter inequality is essential when applying \cite{doucmoulines2011}, Proposition 5; see also the proof of \cite{doucmoulines2011}, Proposition 8.
\end{remark}
%
%
\begin{remark}
In the case $r = 1$ we may replace \hyperref[assum:likelihoodDrift]{(A1)} by the
simpler assumption
that there exists a set $\set{K} \in\field{Y}$ such that the following holds:
\begin{longlist}[(iii)]
\item[(i)]
$
\prob{Y_0 \in\set{K}} > 2 / 3$.
\item[(ii)] For all $\eta> 0$ there exists a local Doeblin set $\set{C} \in
\field{X}$ such that for all $y \in\set{K}$,
%
\begin{equation}
\label{eqbound-eta-G} \sup_{x \in\set{C}^c} g(x, y) \leq\eta\supnormantri{\llhd{y}}
< \infty.
\end{equation}
\item[(iii)]\label{itemmino-g-simple} There exists a set $\set{D} \in
\alg{X}$ satisfying
\[
\inf_{x \in\set{D}} \Q(x,\set{D}) > 0 \quad\mbox{and}\quad\espBig{\ln ^-
\inf_{x \in\set{D}} g(x, Y_0)} < \infty.  %
\]
\end{longlist}
\end{remark}

For the integer $r \in\nsetpos$ and the set $\set{D} \in\field{X}$
given in \hyperref[assum:likelihoodDrift]{(A1)}, define
%
\begin{eqnarray}\label{eqmort-de-rire}
&& \mdr
\nonumber\\[-9pt]\\[-9pt]
&&\qquad \eqdef \bigl\{ \chi\in\probmeas{\field{X}} \dvtx  \espbig{\ln^- \chi
\Pblockantri{\chunk{Y} {0} {\ell- 1}} \1_{\set{D}}} < \infty\mbox{ for all }
\ell\in\{0, \ldots, r\} \bigr\}.\nonumber
\end{eqnarray}

A simple sufficient condition can be proposed to ensure that $\chi\in
\mathcal{M}(\set{D}, r)$.

%
\begin{proposition}
\label{propsufficient-condition-mdr}
Assume that there exist sets $\set{D}_u \in\field{X}$, $u \in\{0,\ldots, r-1\}$, such that (setting $\set{D}_r= \set{D}$ for notational
convenience) for some $\delta> 0$,
%
\begin{equation}
\label{eqCS2-mort-de-rire} \inf_{x \in\set{D}_{u-1}} \Q(x, \set{D}_u) \geq
\delta,  \qquad u \in \{1, \ldots, r\}
\end{equation}
and
%
\begin{equation}
\label{eqCS1-mort-de-rire} \espBig{\ln^- \inf_{x \in\set{D}_u} g(x,Y_0)} <
\infty,  \qquad u \in\{ 0,\ldots,r\}.
\end{equation}
Then any initial distribution $\chi\in\probmeas{\field{X}}$
satisfying $\chi(\set{D}_0) > 0$ belongs to $\mathcal{M}(\set{D},r)$.
\end{proposition}

%
\begin{remark}
\label{remcheck-CS1-mort-de-rire}
To check \eqref{eqCS1-mort-de-rire} we typically assume that for any
given $y \in\set{Y}$,
$\llhd{y}$ is continuous and that the sets $\set{D}_i$, $i \in\{0,\ldots,r-1\}$, are compact.
This condition then translates into an assumption on some generalized
moments of the process $(Y_n)_{n \in\zset}$.
\end{remark}

%
\begin{remark}
\label{remcheck-CS2-mort-de-rire}
Assume that $\set{X}= \rset^d$ for some $d \in\nsetpos$ (or more
generally, $\set{X}$ is a locally compact separable metric space) and
that $\field{X}$ is the associated Borel \mbox{$\sigma$-}field. Assume in
addition that for any open set $\set{O} \in\field{X}$, the function $x
\to Q(x,\set{O})$ is lower semi-continuous on the space $\set{X}$. Then
for any $\delta> 0$ and any compact set $\set{D}_0 \in\field{X}$,
there exist compact sets $\set{D}_u$, $u \in\{0,\ldots,r-1\}$,
satisfying \eqref{eqCS2-mort-de-rire}.
\end{remark}

We are now ready to state our main results.

\subsection{Tightness of asymptotic variance}


%
\begin{theorem} \label{thmtightness-particle-filter}
Assume \textup{\hyperref[assum:likelihoodDrift]{(A1)}--\hyperref[assum:majo-g]{(A2)}}. Then for all $\chi\in
\mdr$ and all $h \in\bmf{\field{X}}$, the sequence $(\variance
[2]{\chunk{Y}{0}{n-1}}{h})_{n \in\nsetpos}$ [defined in~\eqref{eqexpression-asymptotic-variance}] is tight.
\end{theorem}
\begin{pf} 
Using definition \eqref{eqexpressionpred} of the predictive
distribution and the decomposition \eqref{eqLkh-composition} of the
likelihood, we get for all $k \in\{0, \ldots, n - 1\}$,
\[
\predantri[\chi]{\chunk{Y} {0} {n-1}}h = \frac{\chi\Lkh{\chunk{Y}{0}{n-1}}h}{\chi\Lkh{\chunk{Y}{0}{n-1}} \1
_{\set{X}}} = \frac{\chi\Lkh{\chunk{Y}{0}{k-1}} \Lkh{\chunk{Y}{k}{n-1}}h}{\chi\Lkh
{\chunk{Y}{0}{k-1}} \Lkh{\chunk{Y}{k}{n-1}} \1_{\set{X}}}.
\]
Plugging this identity into the expression \eqref
{eqexpression-asymptotic-variance} of the asymptotic variance yields
%
\begin{equation}
\label{eqvariance}\quad  \varianceantri[2]{\chunk{Y} {0} {n-1}} {h} = \sum
_{k = 0}^n \int\predantri[\chi]{\chunk{Y} {0} {k-1}}(\mathrm{d}
x) \biggl[ \frac{\DDelta{\delta_x,\pred[\chi]{\chunk{Y}{0}{k-1}}}{\chunk{Y}{k}{n -
1}}{h,\1_{\set{X}}}} {
(\pred[\chi]{\chunk{Y}{0}{k-1}} \Lkh{\chunk{Y}{k}{n-1}} \1_{\set{X}}
)^2} \biggr]^2,
\end{equation}
where for all sequences $\chunk{y}{k}{n-1} \in\set{Y}^{n - k}$,
functions $f$ and $h$ in $\bmf{\field{X}}$ and probability
measures $\chi$ and $\chi'$ in $\probmeas{\field{X}}$,
%
\begin{eqnarray}\label{eqdef-Delta}
&& \DDeltaantri{\chi,\chi'} {\chunk{y} {k} {n-1}} {f,h}
\nonumber\\[-8pt]\\[-8pt]
&&\qquad \eqdef\chi\Lkhantri{\chunk{y} {k} {n-1}}f \times\chi' \Lkhantri{\chunk{y} {k}
{n-1}}h
- \chi\Lkhantri{\chunk{y} {k} {n-1}}h \times\chi' \Lkhantri{\chunk {y} {k}
{n-1}}f.\nonumber
\end{eqnarray}
Using \eqref{eqexpressionpred}, we obtain for all sequences $\chunk
{y}{0}{n-1} \in\set{Y}^n$,
\begin{eqnarray*}
\predantri[\chi]{\chunk{y} {0} {k - 1}} \Lkhantri{\chunk{y} {k} {n - 1}}
\1_{\set{X}} &=& \frac{\chi\Lkh{\chunk{y}{0}{n-1}}\1_{\set{X}}}{\chi\Lkh{\chunk
{y}{0}{k-1}}\1_{\set{X}}}
\\
&=& \prod_{\ell= k}^{n - 1} \frac{\chi\Lkh{\chunk{y}{0}{\ell}}\1_{\set
{X}}}{\chi\Lkh{\chunk{y}{0}{\ell-1}} \1_{\set{X}}} =
\prod_{\ell
=k}^{n-1} \likeliantri[\chi]{\chunk{y} {0}{\ell-1}} {y_\ell},
\end{eqnarray*}
where $\likeli[\chi]{\chunk{y}{0}{\ell-1}}{y_\ell}$ is the density of
the conditional distribution of $Y_\ell$ given $\chunk{Y}{0}{\ell- 1}$
(i.e., the one-step observation predictor at time $\ell$) defined by
%
\begin{equation}
\label{eqdefinition-likelihood} \likeliantri[\chi]{\chunk{y} {0} {\ell-1}} {y_\ell} \eqdef
\int\predantri[\chi ]{\chunk{y} {0} {\ell-1}}(\mathrm{d} x) g(x,y_{\ell}).
\end{equation}
With this notation, the likelihood function $\chi\Lkh{\chunk{y}{0}{n -
1}} \1_{\set{X}}$ equals the product $\prod_{k = 0}^{n - 1} \likeli[\chi
]{\chunk{y}{0}{k - 1}}{y_k}$ [where we let $\likeli[\chi]{\chunk
{y}{0}{- 1}}{y_0}$ denote the marginal density of~$Y_0$].\vspace*{1pt}

Now, using coupling results obtained in \cite{doucmoulines2011} one
may prove that the predictor distribution forgets its initial
distribution exponentially fast under the $r$-local Doeblin assumption
\eqref{eqdefinition-LD-set}. Moreover, this implies that also the
log-density of the one-step observation predictor forgets its initial
distribution exponentially fast; that is, for all initial distributions
$\chi$ and $\chi'$ there is a deterministic constant $\beta\in\,\ooint
{0,1}$ and an almost surely bounded random variable $C_{\chi,\chi'}$
such that for all $(k, m) \in\nsetpos\times\nset$ and almost all
observation sequences,
%
\begin{equation}
\label{eq1-steppredmixing} \bigl\llvert \ln\likeliantri[\chi]{\chunk{Y} {-m} {k - 1}}
{Y_k} - \ln\likeliantri[\chi ']{\chunk{Y} {-m} {k
- 1}} {Y_k} \bigr\rrvert \leq C_{\chi, \chi'} \beta^{k +
m}.
\end{equation}
Using this, it is shown in \cite{doucmoulines2011}, Proposition~1, that:
\begin{longlist}[(ii)]
\item[(i)]\label{itemlim-ponct-cond-rate-entropy-Y} There exists a
function $\pi\dvtx  \set{Y}^{\zset^{-}} \times\set{Y} \to\rset$ such that
for all probability measures $\chi\in\mdr$,
\[
\lim_{m \to\infty} \likeliantri[\chi]{\chunk{Y} {-m} {-1}}
{Y_0} = \likeliantri{\chunk{Y} {-\infty} {-1}} {Y_0},
\qquad\probalone\mbox{-a.s.} %
\]
Moreover,
%
\begin{equation}
\label{eqesperance-log-stat} \espbig{\bigl|\ln\likeliantri{\chunk{Y} {-\infty} {-1}} {Y_0}\bigr|} <
\infty.
\end{equation}
\item[(ii)] For all probability measures $\chi\in\mdr$, the normalized
log-likelihood function converges according to
%
\begin{equation}
\label{itemlim-ponct-rate-entropy-Y} \lim_{n \to\infty} n^{-1} \ln\chi\Lkhantri{
\chunk{Y} {0} {n-1}} \1_{\set
{X}} = \ell_\infty,  \qquad\probalone\mbox{-a.s.},
\end{equation}
where $\ell_\infty$ is the negated relative entropy, that is, the
expectation of\break  $\ln\likeli{\chunk{Y}{- \infty}{- 1}}{Y_0}$ under the
stationary distribution, that is,
%
\begin{equation}
\label{eqdefinitionellinfty} \ell_\infty\eqdef\espbig{\ln\likeliantri{\chunk{Y} {-\infty}
{-1}} {Y_0}}.
\end{equation}
\end{longlist}
As a first step, we bound the asymptotic variance $\variance
[2]{h}{\chunk{Y}{0}{n - 1}}$ [defined in~\eqref
{eqexpression-asymptotic-variance}] by the product of two quantities,
namely $\variance[2]{\chunk{Y}{0}{n-1}}{h} \leq A \times B_n$, where
%
\begin{eqnarray}
\label{eqdefinition-A} A &\eqdef& \Biggl(\sup_{(k,m) \in\nset^2 \dvtx    k \leq m} \prod
_{\ell=
k}^{m - 1} \frac{\likeli{\chunk{Y}{-\infty}{\ell-1}}{Y_\ell}}{
\likeli[\chi]{\chunk{Y}{0}{\ell- 1}}{Y_\ell}} \Biggr)^4,
\\
\label{eqdefinition-B} B_n &\eqdef& \sum_{m = 0}^n
\biggl( \frac{\sup_{x \in\set{X}} |\DDelta{\delta_x,\predantri[\chi]{\chunk
{Y}{0}{m-1}}}{\chunk{Y}{m}{n-1}}{h,\1_\set{X}}|}{[\prod_{\ell=m}^{n-1}
\likeli{\chunk{Y}{-\infty}{\ell-1}}{Y_\ell}]^2} \biggr)^2.
\end{eqnarray}
Quantity \eqref{eqdefinition-A} can be bounded using the exponential
forgetting \eqref{eq1-steppredmixing} of the one-step predictor
log-density. More precisely, note that
\[
\likeliantri[\chi]{\chunk{Y} {-m} {\ell-1}} {Y_\ell}=\frac{\chi\Lkh{\chunk
{Y}{-m}{\ell}}\1_\set{X}}{\chi\Lkh{\chunk{Y}{-m}{\ell-1}}\1_\set{X}};
\]
thus, by applying Proposition~\ref{propbound-difference-improved}(ii), we conclude that
there exist $\beta\in\,\ooint{0,1}$ and a $\probalone$\mbox{-a.s.} finite
random variable $C_{\chi}$ such that for all $n \in\nset$,
%
\begin{eqnarray}
\label{eqAbound}
&& \prod_{\ell= k}^n
\frac{\likeli{\chunk{Y}{-\infty}{\ell-1}}{Y_\ell
}}{\likeli[\chi]{\chunk{Y}{0}{\ell- 1}}{Y_\ell}}\nonumber
\\
&&\qquad = \prod_{\ell= k}^n \prod
_{m = 0}^\infty\frac{\likeli[\chi]{\chunk{Y}{-
m - 1}{\ell- 1}}{Y_\ell}}{\likeli[\chi]{\chunk{Y}{-m}{\ell- 1}}{Y_\ell}}
\leq\prod_{\ell= k}^n \prod
_{m=0}^\infty\exp\bigl(C_\chi
\beta^{\ell+
m}\bigr)
\\
&&\qquad \leq \exp\bigl(C_\chi/(1 -
\beta)^2\bigr) < \infty,  \qquad\probalone\mbox{-a.s.},\nonumber
\end{eqnarray}
implying that $A$ is indeed $\probalone\mbox{-a.s.}$ finite.

Consider now the second quantity \eqref{eqdefinition-B}. Since the
process $\sequence{Y}[n][\zset]$ is strictly stationary, $\chunk
{Y}{0}{n - 1}$ has the same distribution as $\chunk{Y}{-n}{-1}$ for all
$n \in\nsetpos$. Therefore, for all $n \in\nsetpos$, the random
variable $B_n$ has the same distribution as
%
\begin{equation}
\label{eqdefinition-tilde-B} \widetilde{B}_n \eqdef\sum
_{m = 0}^n \biggl( \frac{\sup_{x \in\set{X}}|\DDelta{\delta_x,\predantri[\chi]{\chunk
{Y}{-n}{-m-1}}}{\chunk{Y}{-m}{-1}}{h,\1_\set{X}}|} {
[\prod_{\ell= 1}^m \likeli{\chunk{Y}{-\infty}{-\ell- 1}}{Y_{-\ell
}}]^2 }
\biggr)^2.
\end{equation}
We will show that $\sup_{n \in\nsetpos} \widetilde{B}_n$ is $\probalone$\mbox{-a.s.}
finite, which implies that the sequence $\sequence{B}[n][\nsetpos]$
is tight. We split each term of $\widetilde{B}_n$ into two factors
according to
%
\begin{eqnarray}\label{eqdecomposition}
&& \frac{\sup_{x \in\set{X}}| \DDelta{\delta_x, \pred[\chi]{\chunk
{Y}{-n}{-m - 1}}}{\chunk{Y}{-m}{-1}}{h, \1_\set{X}}|} {
[\prod_{\ell= 1}^m \likeli{\chunk{Y}{-\infty}{-\ell- 1}}{Y_{-\ell
}}]^2}
\nonumber\\[-8pt]\\[-8pt]
&&\qquad = \biggl(\frac{\supnormTxt{\Lkh{\chunk{Y}{-m}{-1}} \1_\set{X}}}{\prod_{\ell= 1}^m \likeli{\chunk{Y}{-\infty}{-\ell- 1}}{Y_{-\ell}}} \biggr)^2 \frac{\sup_{x \in\set{X}}| \DDelta{\delta_x, \pred[\chi]{\chunk
{Y}{-n}{-m - 1}}}{\chunk{Y}{-m}{-1}}{h, \1_\set{X}}|} {
\| \Lkh{\chunk{Y}{-m}{-1}} \1_{\set{X}} \|^2_{\infty} }\nonumber
\end{eqnarray}
and consider each factor separately.

We will show that the first factor in \eqref{eqdecomposition} grows at
most subgeometrically fast. Indeed, note that
\[
\label{eqfirst-term} \biggl( \frac{\supnormTxt{\Lkh{\chunk{Y}{-m}{-1}} \1_\set{X}}}{\prod_{\ell= 1}^m \likeli{\chunk{Y}{-\infty}{-\ell- 1}}{Y_{-\ell}}} \biggr)^2 = \exp(m
\varepsilon_m),
\]
where
\[
\label{eqdefinition-varepsilon} \varepsilon_m \eqdef\frac{2}{m} \Biggl( \ln
\supnormantri{\Lkhantri{\chunk {Y} {-m} {-1}} \1_\set{X}} - \sum
_{\ell= 1}^m \ln\likeliantri{\chunk{Y} {-\infty} {-\ell-1}}
{Y_{-\ell
}} \Biggr).
\]
According to Lemma~\ref{lemmekong}, $\varepsilon_m \to2(\ell
_\infty- \ell_\infty) = 0$, \probalone\mbox{-a.s.}, as $m \rightarrow\infty$.

The second factor in \eqref{eqdecomposition} is handled using
Proposition~\ref{propbound-difference-improved}(iii),
which guarantees the existence of a constant $\beta\in\,\ooint{0,1}$
and a $\probalone$\mbox{-a.s.} random variable $C$ such that for all $(m,n) \in
(\nsetpos)^2$,
%
\begin{equation}
\label{eqBbound} \frac{\sup_{x \in\set{X}}| \DDelta{\delta_x, \pred[\chi]{\chunk
{Y}{-n}{-m - 1}}}{\chunk{Y}{-m}{-1}}{h, \1_\set{X}}|}{\| \Lkh{\chunk
{Y}{-m}{-1}} \1_{\set{X}} \|^2_{\infty} } \leq C \beta^m \supnorm{h}.
\end{equation}
This concludes the proof.
\end{pf}


%
\begin{theorem} \label{thmtightness-posterior}
Assume \textup{\hyperref[assum:likelihoodDrift]{(A1)}--\hyperref[assum:majo-g]{(A2)}}. Then for all $\chi\in
\mdr$ and all $h \in\bmf{\field{X}}$, the sequence $(\filtvariance
[2]{\chunk{Y}{0}{n}}{h})_{n \in\nsetpos}$ [defined in~\eqref
{eqdeffiltvariance}] is tight.
\end{theorem}

\begin{pf}
Using the expression \eqref{eqvariance} of the asymptotic variance of
the predictor approximation yields for all $\chunk{y}{0}{n} \in\set
{Y}^{n + 1}$, as
\[
\predantri[\chi]{\chunk{y} {0} {n - 1}} \llhd{y_n} = \likeliantri[\chi]{\chunk
{y} {0} {n - 1}} {y_n} %
\]
and
\[
\likeliantri[\chi]{\chunk{y} {0} {n - 1}} {y_n} \times\predantri[\chi]{\chunk
{y} {0} {k - 1}} \Lkhantri{\chunk{y} {k} {n - 1}} \1_{\set{X}} = \predantri[\chi ]{
\chunk{y} {0} {k - 1}} \Lkhantri{\chunk{y} {k} {n}} \1_{\set{X}},  %
\]
the identity
\begin{eqnarray*}
&& \filtvarianceantri[2]{\chunk{y} {0} {n}} {h}
\\
&&\qquad = \sum_{k = 0}^n \int\predantri[\chi ]{
\chunk{y} {0} {k - 1}}(\mathrm{d} x) \biggl[ \frac{\DDelta{\delta_x, \pred[\chi
]{\chunk{y}{0}{k - 1}}}{\chunk{y}{k}{n - 1}}{\llhd{y_n} \{ h - \filt
[\chi]{\chunk{y}{0}{n}}h \}, \1_{\set{X}}}}{\pred[\chi]{\chunk{y}{0}{k
- 1}} \Lkh{\chunk{y}{k}{n - 1}} \1_{\set{X}} \times\pred[\chi]{\chunk
{y}{0}{k - 1}} \Lkh{\chunk{y}{k}{n}} \1_{\set{X}}} \biggr]^2,
\end{eqnarray*}
where
%
\begin{eqnarray}
\label{eqdiff} && \DDeltaantri{\delta_x, \pred[\chi]{\chunk{y} {0} {k -1}}}
{\chunk{y} {k} {n - 1}} {\llhd{y_n} \bigl\{ h - \filtantri[\chi]{\chunk{y} {0} {n}}h \bigr\}, \1_{\set{X}}}
\nonumber
\\
&&\qquad= \delta_x \Lkhantri{\chunk{y} {k} {n - 1}} \bigl(
\llhd{y_n} \bigl\{ h - \filtantri[\chi]{\chunk{y} {0} {n}}h \bigr\}\bigr)
\times\predantri[\chi]{\chunk{y} {0} {k - 1}} \Lkhantri{\chunk{y} {k} {n - 1}}
\1_{\set{X}}
\\
&&\quad\qquad{} - \delta_x \Lkhantri{\chunk{y} {k} {n - 1}}
\1_{\set{X}} \times\predantri[\chi]{\chunk{y} {0} {k - 1}} \Lkhantri{\chunk{y} {k} {n
- 1}} \bigl(\llhd {y_n} \bigl\{ h - \filtantri[\chi]{\chunk{y} {0} {n}}h
\bigr\}\bigr).
\nonumber
\end{eqnarray}
Here the equation
%
\begin{equation}
\label{eqfiltidentity} \filtantri[\chi]{\chunk{y} {0} {n}}h = \frac{\pred[\chi]{\chunk{y}{0}{k - 1}}
\Lkh{\chunk{y}{k}{n - 1}}(\llhd{y_n} h)}{\pred[\chi]{\chunk{y}{0}{k -
1}} \Lkh{\chunk{y}{k}{n - 1}}\llhd{y_n}}
\end{equation}
implies that
\begin{eqnarray*}
&& \predantri[\chi]{\chunk{y} {0} {k - 1}} \Lkhantri{\chunk{y} {k} {n - 1}} \bigl(
\llhd{y_n} \bigl\{ h - \filtantri[\chi]{\chunk{y} {0} {n}}h \bigr\}\bigr)
\\
&&\qquad = \predantri[\chi]{\chunk{y} {0} {k - 1}} \Lkhantri{\chunk{y} {k} {n - 1}} \bigl(
\llhd{y_n} h\bigr)
- \filtantri[\chi]{\chunk{y} {0} {n}}h \times\predantri[\chi]{\chunk{y} {0} {k - 1}}
\Lkhantri{\chunk{y} {k} {n - 1}} \llhd{y_n}
\\
&&\qquad = 0,
\end{eqnarray*}
which implies in turn that the second term on the right-hand side of
\eqref{eqdiff} vanishes. Thus, developing also the first term and
reusing the identity \eqref{eqfiltidentity} yields
\begin{eqnarray*}
&& \frac{\DDelta{\delta_x, \pred[\chi]{\chunk{y}{0}{k - 1}}}{\chunk
{y}{k}{n - 1}}{\llhd{y_n} \{ h - \filt[\chi]{\chunk{y}{0}{n}}h \}, \1
_{\set{X}}}}{\pred[\chi]{\chunk{y}{0}{k - 1}} \Lkh{\chunk{y}{k}{n - 1}}
\1_{\set{X}} \times\pred[\chi]{\chunk{y}{0}{k - 1}} \Lkh{\chunk
{y}{k}{n}} \1_{\set{X}}}
\\
&&\qquad= \frac{\delta_x \Lkh{\chunk{y}{k}{n - 1}} (\llhd{y_n} h) -
\filt[\chi]{\chunk{y}{0}{n}}h \times\delta_x \Lkh{\chunk{y}{k}{n - 1}}
\llhd{y_n}}{\pred[\chi]{\chunk{y}{0}{k - 1}} \Lkh{\chunk{y}{k}{n}} \1
_{\set{X}}}
\\
&&\qquad= \frac{\DDelta{\delta_x, \pred[\chi]{\chunk{y}{0}{k -
1}}}{\chunk{y}{k}{n - 1}}{\llhd{y_n} h, \llhd{y_n}}}{\pred[\chi]{\chunk
{y}{0}{k - 1}} \Lkh{\chunk{y}{k}{n - 1}}\llhd{y_n} \times\pred[\chi
]{\chunk{y}{0}{k - 1}} \Lkh{\chunk{y}{k}{n}} \1_{\set{X}}}
\\
&&\qquad= \frac{\DDelta{\delta_x, \pred[\chi]{\chunk{y}{0}{k -
1}}}{\chunk{y}{k}{n - 1}}{\llhd{y_n} h, \llhd{y_n}}}{(\pred[\chi]{\chunk
{y}{0}{k - 1}} \Lkh{\chunk{y}{k}{n}} \1_{\set{X}})^2}.
\end{eqnarray*}
Thus, to sum up,
\[
\filtvarianceantri[2]{\chunk{y} {0} {n}} {h}
\\
= \sum_{k = 0}^n \int\predantri[\chi ]{
\chunk{y} {0} {k - 1}}(\mathrm{d} x) \biggl[ \frac{\DDelta{\delta_x, \pred[\chi
]{\chunk{y}{0}{k - 1}}}{\chunk{y}{k}{n - 1}}{\llhd{y_n} h, \llhd
{y_n}}}{(\pred[\chi]{\chunk{y}{0}{k - 1}} \Lkh{\chunk{y}{k}{n}} \1_{\set
{X}})^2} \biggr]^2,
\]
providing an expression of the asymptotic variance of the particle
filter that resembles closely the corresponding expression \eqref
{eqvariance} for the particle predictor. Using this, the tightness can
be established along the very same lines as Theorem~\ref{thmtightness-particle-filter}, and we leave the details to the
interested reader.
\end{pf}

\subsection{Tightness of the asymptotic $\mathsf{L}^{p}$ error}

In the following we show that tightness of the asymptotic variance
implies tightness of the asymptotic $\Lp{p}$ error (when scaled with
$\sqrt{N}$). The asymptotic $\Lp{p}$ error given in Theorem~\ref{thmasymptoticLpbound} below is obtained by establishing, for fixed
time indices $n$, using a standard exponential deviation inequality,
uniform integrability (with respect to the particle sample size $N$) of
the sequence of scaled $\Lp{p}$ errors. After this, weak convergence
implies convergence of moments, implying in turn convergence of the $\Lp
{p}$ error.

%
\begin{theorem} \label{thmasymptoticLpbound}
Assume \textup{\hyperref[assum:majo-g]{(A2)}}. Then, for all functions $h \in\bmf{\field
{X}}$, constants \mbox{$p \in\rsetpos$} and all initial distributions $\chi
\in\probmeas{
\alg{X}}$, $\probalone\mbox{-a.s.}$,
%
\begin{eqnarray*}
&& \lim_{N \to\infty} \sqrt{N} \cespbig{\bigl\llvert \partpredantri[\chi]{\chunk
{Y} {0} {n - 1}} h - \predantri[\chi]{\chunk{Y} {0} {n-1}} h \bigr\rrvert
^p} {\chunk {Y} {0} {n-1}}[1/p]
\\
&&\qquad = \sqrt{2} \varianceantri{\chunk{Y} {0} {n-1}} {h} \biggl( \frac{\Gamma
((p+1)/2)}{\sqrt{2 \pi}}
\biggr)^{1/p},
\end{eqnarray*}
where $\Gamma$ is the gamma function.
\end{theorem}

\begin{pf}
Recall that if $\sequence{A}[N][\nsetpos]$ is a sequence of random
variables such that $A_N \dlim A$ as $N \to\infty$ and $\sequence
{A^p}[N][\nsetpos]$ is uniformly integrable for some $p> 0$, then
$\espTxt{|A|^p} < \infty$, $\lim_{N \to\infty} \espTxt[]{A_N^p} =
\espTxt[]{A^p}$ and $\lim_{N \to\infty} \espTxt[]{|A_N|^p}= \espTxt
[]{|A|^p}$; see, for example, \cite{serfling1980}, Theorem~A, page~14.
Now set, for $n \in\nsetpos$,
\[
A_{N, \chi} \bigl\langle\chunk{Y} {0} {n - 1} \bigr\rangle(h) \eqdef\sqrt{N} \bigl(
\partpredantri[\chi]{\chunk{Y} {0} {n - 1}}h - \predantri[\chi]{\chunk{Y} {0} {n - 1}}h
\bigr).  %
\]
Let $q > p$ and write
\begin{eqnarray*}
&& \sup_{N \in\nsetpos} \cespbig{\bigl\llvert A_{N,\chi} \bigl\langle
\chunk {Y} {0} {n-1}\bigr\rangle(h)\bigr\rrvert ^q} {\chunk{Y} {0} {n-1}}
\\
&&\qquad= \sup_{N \in\nsetpos} \int_0^\infty
\probantri{\bigl\llvert A_{N,\chi} \bigl\langle\chunk{Y} {0} {n-1}\bigr\rangle(h) \bigr
\rrvert \geq \epsilon^{1/q}}[\chunk{Y} {0} {n-1}]\,\mathrm{d}\epsilon
\\
&&\qquad= q \sup_{N \in\nsetpos} \int_0^\infty
\epsilon^{q-1}\probantri{\bigl\llvert A_{N,\chi} \bigl\langle
\chunk{Y} {0} {n-1}\bigr\rangle(h)\bigr\rrvert \geq \epsilon}[\chunk{Y} {0} {n-1}]\,\mathrm{d}\epsilon.
\end{eqnarray*}
Now, note that \hyperref[assum:majo-g]{(A2)}(ii) implies that
$\supnormTxt{\llhd{Y_n}}$ is $\probalone$\mbox{-a.s.} finite for all $n \in
\nset$. Thus, Assumptions~1 and 2 of \cite{doucgariviermoulinesolsson2009} are fulfilled, which implies, via
Remark~1 in the same work (see also \cite{doucguillinnajim2005}, Lemma 2.1, \cite{delmoraljacodprotter2001}, Theorem~3.1, \cite{delmoralmiclo2000}, Theorem 3.39, and \cite{delmoralledoux2000}, Lemma~4, for
similar results), that there exist, for all $n \in\nset$ and $h \in
\bmf{\field{X}}$, positive constants $B_n$ and $C_n$ (where only the
latter depends on $h$) such that for all $N \in\nsetpos$ and all
$\epsilon> 0$,
%
\begin{equation}
\label{eqhoeffdinginequalitypredictor}\probantri{\bigl|A_{N, \chi} \bigl\langle\chunk{Y} {0} {n - 1}
\bigr\rangle(h)\bigr| \geq\epsilon }[\chunk{Y} {0} {n - 1}] \leq B_n \exp
\bigl(-C_n \epsilon^2\bigr).
\end{equation}
This implies that for all $n \in\nset$, \probalone\mbox{-a.s.},
\[
\sup_{N \in\nsetpos} \cespbig{\bigl\llvert A_{N,\chi} \bigl\langle\chunk
{Y} {0} {n-1}\bigr\rangle(h)\bigr\rrvert ^q} {\chunk{Y} {0} {n-1}} \leq q
B_n \int_0^\infty\epsilon^{q-1}
\exp\bigl(-C_n \epsilon^2\bigr)\,\mathrm{d}\epsilon< \infty,
\]
which establishes, via \cite{shiryaev1996}, Lemma~II.6.3, page~190, as
$q > p$, that the sequence $( |A_{N, \chi} \langle\chunk{Y}{0}{n - 1}
\rangle(h)|^p)_{N \in\nset}$ is uniformly integrable conditionally on
$\chunk{Y}{0}{n-1}$, that is,
\[
\lim_{M \to\infty} \sup_{N \in\nsetpos} \cespbig{\bigl\llvert
A_{N,\chi} \bigl\langle\chunk{Y} {0} {n - 1} \bigr\rangle(h)\bigr\rrvert
^p \1_{\{| A_{N, \chi}
\langle\chunk{Y}{0}{n - 1}\rangle(h) | \geq M \}}} {\chunk{Y} {0} {n-1}} = 0,  \qquad\probalone
\mbox{-a.s.}
\]
We may now complete the proof by applying Theorem~\ref{thmcltparticle}, which states that conditionally on $\chunk{Y}{0}{n
- 1}$, as $N \rightarrow\infty$,
\[
A_{N, \chi} \bigl\langle\chunk{Y} {0} {n - 1} \bigr\rangle(h) \dlim\varianceantri {
\chunk{Y} {0} {n - 1}} {h} Z,
\]
where $Z$ is a standard normally distributed random variable.
\end{pf}

We next state the corresponding result for the particle filter
approximation, which is obtained along the very same lines as
Theorem~\ref{thmasymptoticLpbound}.

%
\begin{theorem} \label{thmasymptoticLpboundposterior}
Assume \textup{\hyperref[assum:majo-g]{(A2)}}. Then, for all functions $h \in\bmf{\field
{X}}$, constants \mbox{$p \in\rsetpos$} and
all initial distributions $\chi\in\probmeas{\field{X}}$, $\probalone
\mbox{-a.s.}$,
%
\begin{eqnarray*}
&& \lim_{N \to\infty} \sqrt{N} \cespbig{\bigl\llvert \partfiltantri[\chi]{\chunk
{Y} {0} {n}} h - \filtantri[\chi]{\chunk{Y} {0} {n}} h \bigr\rrvert ^p} {
\chunk {Y} {0} {n-1}}[1/p]
\\
&&\qquad = \sqrt{2} \filtvarianceantri{\chunk{Y} {0} {n-1}} {h} \biggl( \frac{\Gamma
((p+1)/2)}{\sqrt{2 \pi}}
\biggr)^{1/p},
\end{eqnarray*}
where $\Gamma$ is the gamma function.
\end{theorem}

\section{Examples}\label{secapplications}
In this section, we develop two classes of examples. In Section~\ref{secpostcompkalman} we consider the \emph{linear Gaussian state--space
models}, an important model class that is used routinely in time-series
analysis. Recall that in the linear Gaussian case, closed-form
solutions to the optimal filtering problem can be obtained using the
Kalman recursions. However, as an illustration, we analyze this model
class under assumptions that are very general. In Section~\ref{secexample-nonlinear-state-space}, we consider a significantly more
general class of nonlinear state--space models. In both these examples
we will find that assumptions~\hyperref[assum:likelihoodDrift]{(A1)}--\hyperref[assum:majo-g]{(A2)}
are satisfied and straightforwardly verified.

\subsection{Linear Gaussian state--space models}\label{secpostcompkalman}

The linear Gaussian state--space models form an important class of HMMs.
Let $\set{X} = \rset^{d_x}$ and $\set{Y} = \rset^{d_y}$ and define
state and observation sequences through the linear dynamic system
\begin{eqnarray*}
X_{k+1} & =& A X_k + \URoot U_k,
\\
Y_k & =& B X_k + \VRoot V_k,
\end{eqnarray*}
where $(U_k, V_k)_{k \geq0}$ is a sequence of i.i.d. Gaussian
vectors with zero mean and identity covariance matrix. The noise
vectors are assumed to be independent of~$X_0$.
Here $U_k$ is $d_u$-dimensional, $V_k$ is $d_y$-dimensional and the
matrices $A$, $R$, $B$ and $S$ have the appropriate dimensions. Note
that we cover also the case $d_u < d_x$, for which the prior kernel $\Q
$ does not admit a transition density with respect to Lebesgue measure.

For any $n \in\nset$, define the observability and controllability
matrices $\mathcal{O}_{n}$ and $\mathcal{C}_{n}$ by
%
\begin{equation}
\label{eqdefinition-controlabilite-observabilite} \mathcal{O}_{n}\eqdef \lleft[\matrix{ B
\cr
B A
\cr
B A^2
\cr
\vdots
\cr
B A^{n-1}} \rright] \quad
\mbox{and}\quad \mathcal{C}_{n}\eqdef \bigl[
\matrix{A^{n-1}\URoot& A^{n-2} \URoot&\cdots&\URoot} \bigr],
\end{equation}
respectively. We assume the following.
\begin{longlist}[(LGSS3)]
\item[(LGSS1)]\label{hypLGSScommandable-observable}
The pair $(A, B)$ is observable, and the
pair $(A, \URoot)$ is controllable, that is, there exists $r \in\nset$
such that the observability
matrix $\mathcal{O}_r$ and the controllability matrix $\mathcal{C}_r$
have full rank.

\item[(LGSS2)]\label{hypLGSSobservation-noise-full-rank}
The measurement noise covariance matrix $\VRoot$ has full rank.

\item[(LGSS3)]\label{hypLGSSvariance-observation} $\esp{\|Y_0\|^2} <
\infty$.
\end{longlist}

We now check assumptions~\hyperref[assum:likelihoodDrift]{(A1)}--\hyperref[assum:majo-g]{(A2)}. The
dimension $d_u$ of the state noise vector $U_k$ is in many situations
smaller than the dimension $d_x$ of the state vector $X_k$ and hence
$\UCov$ may be rank deficient (here $\transp{}$ denotes the transpose).
Some additional notation is required: for any positive matrix $A$ and
vector $z$ of appropriate dimension, denote $\| z \|_A^2 \eqdef\transp
{z} A^{-1} z$. In addition, define for any $n \in\nset$,
%
\begin{equation}
\label{eqdefinition-f} \mathcal{F}_{n} \eqdef\mathcal{D}_{n}
\transp{\mathcal{D}_{n}} + \mathcal{S}_{n} \transp{
\mathcal{S}_{n}},
\end{equation}
where
\[
\mathcal{D}_n \eqdef \lleft[\matrix{ 0 & 0 & \cdots & 0
\cr
B R
& \ddots & & 0
\cr
B A R & B R & \ddots & \vdots
\cr
\vdots & & \ddots & 0
\cr
B
A^{n-2} R & B A^{n-3} R & \cdots& B R } \rright], \qquad
\mathcal{S}_n \eqdef \lleft[ \matrix{ S & 0 & \cdots& 0
\cr
0 & S &
\ddots& \vdots
\cr
\vdots& \ddots& \ddots& 0
\cr
0 & \cdots& 0& S} \rright].
\]
Under (LGSS2), the matrix $\mathcal
{F}_{n}$ is positive definite for any $n \geq r$. When the state
process is initialized at $x_0 \in\set{X}$, the likelihood of the
observations $\chunk{y}{0}{n-1} \in\set{Y}^n$ is given by
\[
\label{eqLGSSlikelihood-chunk-observation}
\delta_{x_0} \Pblockantri{\chunk{y} {0} {n-1}}
\1_{\set{X}} = (2 \pi)^{-n d_y} \operatorname{det}^{-1/2}(
\mathcal{F}_{n}) \exp \bigl( - \tfrac{1}2 \|
\mathbf{y}_{n-1} - \mathcal{O}_{n} x_0 \|
^2_{\mathcal{F}_{n}} \bigr),
\]
where $\mathbf{y}_{n-1} \eqdef\transp{[\transp{y_0}, \transp{y_1},\ldots, \transp{y_{n-1}}]}$ and $\mathcal{O}_{n}$ is defined in \eqref
{eqdefinition-controlabilite-observabilite}.

We first consider \hyperref[assum:likelihoodDrift]{(A1)}. Under (LGSS1), the observability matrix $\mathcal
{O}_{r}$ is full rank, and we have for any compact subset $\set{K}
\subset\set{Y}^{r}$,
\[
\lim_{\| x_0 \| \to\infty} \inf_{\chunk{y}{0}{r-1} \in\set{K}} \llVert
\mathbf{y}_{r-1} - \mathcal{O}_{r} x_0 \rrVert
_{\mathcal{F}_{r}} = \infty,
\]
showing that for all $\eta> 0$, we may choose a compact set $\set{C}
\subset\rset^{d_x}$ such that \eqref{eqbound-eta-G22} is satisfied. It
remains to prove that any compact set $\set{C}$ is an $r$-local Doeblin
set satisfying condition \eqref{eqlower-bound}. For any $\chunk
{y}{0}{r - 1} \in\set{Y}^r$ and $x_0 \in\set{X}$, the measure $\delta
_{x_0} \Pblock{\chunk{y}{0}{r - 1}}$\vadjust{\goodbreak} is absolutely continuous
with respect to the Lebesgue measure on $(\set{X}, \field{X})$ with
Radon--Nikodym derivative $\pblock{\chunk{y}{0}{r - 1}}(x_0, x_r)$
given (up to an irrelevant multiplicative factor) by
%
\begin{equation}
\label{eqGLSSradon-nikodym} \pblockantri{\chunk{y} {0} {r-1}}(x_0,x_r)
\propto\operatorname{det}^{-1/2}( \mathcal{G}_{r}) \exp\lleft(
- \frac{1}{2} \biggl\llVert \lleft[ \matrix{ \mathbf{y}_{r-1}
\cr
x_r } \rright] - \lleft[ \matrix{ \mathcal{O}_{r}
\cr
A^r } \rright] x_0\biggr\rrVert _{\mathcal{G}_{r}}^2
\rright),
\end{equation}
where the covariance matrix $\mathcal{G}_r$ is
\[
\mathcal{G}_r \eqdef \lleft[\matrix{ \mathcal{D}_{r}
\cr
\mathcal{C}_{r}} \rright] \bigl[ \matrix{ \transp{
\mathcal{D}_{r}} & \transp{\mathcal{C}_{r}}} \bigr] +
\lleft[ \matrix{ \mathcal{S}_{r}
\cr
\mathbf{0}} \rright] \bigl[
\matrix{ \transp{\mathcal{S}_{r}} &\transp{\mathbf{0}}} \bigr].
\]
The proof of \eqref{eqGLSSradon-nikodym} relies on the positivity of
$\mathcal{G}_{r}$, which requires further discussion.
By construction, the matrix $\mathcal{G}_{r}$ is nonnegative.
For all $\mathbf{y}_{r-1} \in\set{Y}^r$ and $x \in\set{X}$, the equation
\[
\bigl[\matrix{\transp{\mathbf{y}_{r-1}} & \transp{x}}\bigr] \mathcal{G}_r
\lleft[ \matrix{ \mathbf{y}_{r-1}
\cr
x } \rright] = \bigl\| \transp{
\mathcal{D}_r} \mathbf{y}_{r-1} + \transp{\mathcal
{C}_r} x \bigr\| ^2 + \bigl\| \transp{
\mathcal{S}_r} \mathbf{y}_{r-1} \bigr\| ^2 = 0
\]
implies that $\| \transp{\mathcal{D}_{r}} \mathbf{y}_{r-1}+ \transp
{\mathcal{C}_{r}} x \|^2 = 0$ and
$\| \transp{\mathcal{S}_{r}} \mathbf{y}_{r-1} \|^2 =0 $. Since the
matrix $\mathcal{S}_{r}$ has full rank, this implies that
$\mathbf{y}_{r-1} = 0$. Since also $\mathcal{C}_r$ has full rank [the
pair $(A, R)$ is controllable],
this implies in turn that $x = 0$. Therefore, the matrix $\mathcal
{G}_r$ is positive definite and
the function
\[
(x_0,x_r) \mapsto\biggl\llVert \lleft[ \matrix{
\mathbf{y}_{r-1}
\cr
x_r } \rright] - \lleft[ \matrix{
\mathcal{O}_{r}
\cr
A^r } \rright] x_0
\biggr\rrVert _{\mathcal{G}_{r}}^2 %
\]
is continuous for all $\mathbf{y}_{r-1}$. It is therefore bounded on
any compact subset of $\set{X}^2$. This implies that every nonempty
compact set $\set{C} \subset\rset^{d_x}$ is an $r$-local Doeblin set,
with $\lambda_{\set{C}}(\cdot) = \lleb(\cdot) / \lleb(\set{C})$ and
\begin{eqnarray*}
\epsilon^-_{\set{C}}\bigl(\chunk{y} {0} {r-1}\bigr)&=& \bigl( \lleb(\set{C})
\bigr)^{-1} \inf_{(x_0,x_r) \in\set{C}^2} \pblockantri{\chunk {y} {0}
{r-1}}(x_0,x_r),
\\
\epsilon^+_{\set{C}}\bigl(\chunk{y} {0} {r-1}\bigr)&=& \bigl( \lleb(\set{C})
\bigr)^{-1} \sup_{(x_0,x_r) \in\set{C}^2} \pblockantri{\chunk {y} {0}
{r-1}}(x_0,x_r).
\end{eqnarray*}
Consequently, condition \eqref{eqlower-bound} is satisfied for any
compact set $\set{K} \subseteq\set{Y}^{r - 1}$. It remains to verify
\hyperref[assum:likelihoodDrift]{(A1)}(iii). Under (LGSS1), the measure $\delta_{x_0} \Pblock
{\chunk{y}{0}{r-1}}$ is absolutely continuous with respect to the Lebesgue
measure $\lleb$; therefore, for any set $\set{D} \subset\rset^{d_x}$,
\[
\inf_{x_0 \in\set{D}} \delta_{x_0} \Pblockantri{\chunk{y} {0}
{r-1}}(\set{D}) \geq \inf_{(x_0, x_r) \in\set{D}^2} \pblockantri{\chunk{y} {0} {r -
1}}(x_0, x_r) \lleb(\set{D}).
\]
Take $\set{D}$ to be any compact set with positive Lebesgue measure. Now,
\begin{eqnarray*}
&& \sup_{(x_0, x_r) \in\set{D}^2} \biggl\llVert \lleft[ \matrix{
\mathbf{y}_{r-1}
\cr
x_r } \rright] - \lleft[ \matrix{
\mathcal{O}_{r}
\cr
A^r } \rright] x_0
\biggr\rrVert _{\mathcal
{G}_{r}}^2
\\
&&\qquad \leq2 \lambda_{\max} (\mathcal{G}_{r} ) \Bigl\{ \| \mathbf
{y}_{r-1} \|^2 + \max_{x \in\set{D}} \|x
\|^2 \bigl[1+\lambda_{\max
} \bigl(\transp{
\mathcal{O}_{r}} \mathcal{O}_{r} + \transp{A^r}
A^r \bigr) \bigr] \Bigr\},
\end{eqnarray*}
where $\lambda_{\max}(A)$ is the largest eigenvalue of $A$. Under (LGSS3), $\esp{\|Y_0\|^2} < \infty$, implying
that \hyperref[assum:likelihoodDrift]{(A1)}(iii) is
satisfied for any compact set.

We now consider \hyperref[assum:majo-g]{(A2)}. Under (LGSS2), $\VRoot$ has full rank, and
taking the
reference measure $\lleb$ as the Lebesgue measure on $\set{Y}$, $g(x,
y)$ is, for each $x \in\set{X}$, a Gaussian density with covariance
matrix $\VCov$. We therefore have
\[
\supnormantri{\llhd{y}} = (2\pi)^{-d_y/2} \operatorname{det}^{-1/2}\bigl(\VCov\bigr) <
\infty %
\]
for all $y \in\set{Y}$, which verifies \hyperref[assum:majo-g]{(A2)}(i)--(ii).

To conclude this discussion, we need to specify more explicitly the set
$\mathcal{M}(\set{D}, r)$ [see \eqref{eqmort-de-rire}]
of possible initial distributions. Using Proposition~\ref{propsufficient-condition-mdr}, we verify sufficient conditions \eqref
{eqCS2-mort-de-rire} and \eqref{eqCS1-mort-de-rire}. To check \eqref
{eqCS2-mort-de-rire}, we use Remark~\ref{remcheck-CS2-mort-de-rire}: for any open subset $\set{O} \subset\rset
^{d_x}$ and $x \in\set{X}$, $\Q(x,\set{O})= \esp{\1_{\set{O}}(A x + R U)}$,
where the expectation is taken with respect to the $d_u$-dimensional standard
normal distribution. Let $(x_n)_{n \in\nsetpos}$ be a sequence in $\set
{X}$ converging to $x$. By using that the function $\1_{\set{O}}$ is lower
semi-continuous we obtain, via Fatou's lemma,
\[
\liminf_{n \to\infty} \Q(x_n,\set{O}) \geq\espBig{\liminf
_{n \to\infty} \1_{\set{O}}(A x_n + R U)} \geq\Q(x,
\set{O}),
\]
showing that the function $x \mapsto\Q(x,\set{O})$ is lower
semi-continuous for any open subset $\set{O}$.

Assumption (LGSS2) implies that for
all $(x, y) \in\set{X} \times\set{Y}$,
\begin{eqnarray*}
\ln g(x,y) &\geq&-\frac{d_y}{2} \ln( 2 \pi) -\frac{1}2 \ln\operatorname{det}^{-1/2}(\VCov)
\\
&&{}- \bigl[\lambda_{\min} \bigl(\VCov \bigr) \bigr]^{-1} \bigl( \|y
\|^2 + \|B x\|^2 \bigr),
\end{eqnarray*}
where $\lambda_{\min}(\VCov)$ is the minimal eigenvalue of $\VCov$.
Therefore \eqref{eqCS1-mort-de-rire} is satisfied under (LGSS3). Consequently, we may apply
Theorems~\ref{thmtightness-particle-filter}~and~\ref{thmtightness-posterior} to establish tightness of the asymptotic
variances of the particle predictor and filter approximations for any
initial distribution $\chi\in\probmeas{\field{X}}$ as soon as the
process $(Y_k)_{k \in\zset}$ is strictly stationary ergodic and $\esp
{\|Y_0\|^2} < \infty$.

\subsection{Nonlinear state--space models}\label{secexample-nonlinear-state-space}

We now turn to a very general class of nonlinear state--space models.
Let $\set{X}=\rset^d$, $\set{Y}=\rset^\ell$ and $\field{X}$ and $\field
{Y}$ be the associated Borel $\sigma$-fields. In the following we
assume that for each $x \in\set{X}$, the probability measure $\Q(x,
\cdot)$ has a density $q(x, \cdot)$ with respect to the Lebesgue measure $\lleb$
on $\rset^d$. For instance, the state sequence $(X_k)_{k \in\nset}$
could be defined through some nonlinear recursion
%
\begin{equation}
\label{eqARCH} X_k = T (X_{k-1}) + \Sigma(X_{k-1})
\zeta_k,
\end{equation}
where $(\zeta_k)_{k \in\nsetpos}$ is an i.i.d. sequence of
$d$-dimensional random vectors with density $\rho_\zeta$ with respect to the
Lebesgue measure $\lleb$ on $\rset^d$. Here $T \dvtx  \rset^d \to\rset^d$
and $\Sigma\dvtx  \rset^d \to\rset^{d \times d}$ are given (measurable)
matrix-valued functions such that $\Sigma(x)$ is full rank for each $x
\in\set{X}$. Models of form \eqref{eqARCH}, typically referred to as
vector \emph{autoregressive conditional heteroscedasticity} (ARCH) \emph
{models}, are often of interest in time series analysis and financial
econometrics. In this context, we let the observations $(Y_k)_{k \in
\nset}$ be generated through a given measurement density $g(x, y)$
(again with respect to the Lebesgue measure).

We now introduce the basic assumptions of this section.
\begin{longlist}[(NL3)]
\item[(NL1)]\label{itemNLpositivity-density}
The function $(x, x') \mapsto q(x, x')$ on $\set{X}^2$ is positive and
continuous. In addition, $\sup_{(x, x') \in\set{X}^2} q(x, x') < \infty$.

\item[(NL2)]\label{itemNLgoes-to-zero-at-infinity}
For any compact subset $\set{K} \subset\set{Y}$,
\[
\lim_{\| x \| \to\infty} \sup_{y \in\set{K}} \frac{g(x,y)}{\supnormTxt
{\llhd{y}}} =
0.
\]

\item[(NL3)]\label{itemNLbound-log+-g}
For all $(x,y) \in\set{X} \times\set{Y}$, $g(x, y) > 0$ and
\[
\espbig{\ln^+ \supnormantri{\llhd{Y_0}}} < \infty.  %
\]

\item[(NL4)]\label{itemNLbound-log--g}
There exists a compact subset $\set{D} \subset\set{Y}$ such that
\[
\espBig{ \ln^- \inf_{x \in\set{D} } g(x,Y_0)} < \infty.
\]
\end{longlist}

Under (NL1), every compact set $\set{C} \subset\set{X}= \rset^d$ with
positive Lebesgue measure is 1-small and therefore local Doeblin with
$\lambda_{\set{C}}(\cdot) = \lleb(\cdot\cap\set{C})/ \lleb(\set{C})$,
$\varphi_{\set{C}} \langle y_0 \rangle= \lleb(\set{C})$ and
\begin{eqnarray*}
\epsilon_{\set{C}}^- &=& \inf_{(x, x') \in\set{C}^2} q\bigl(x,
x'\bigr),
\\
\epsilon_{\set{C}}^+ &=& \sup_{(x, x') \in\set{C}^2} q\bigl(x,
x'\bigr).
\end{eqnarray*}
Under (NL1) and (NL2), conditions \eqref{eqbound-eta-G} and \eqref
{eqlower-bound} are satisfied with $r = 1$. In addition, \eqref
{eqcondition-minoration} is implied by (NL1)~and~(NL4). Consequently,
assumption \hyperref[assum:likelihoodDrift]{(A1)} holds. Moreover, \hyperref[assum:majo-g]{(A2)}
follows directly from (NL3). So, finally, under (NL1)--(NL4) we
conclude, using Proposition~\ref{propsufficient-condition-mdr},
Theorems~\ref{thmtightness-particle-filter}~and~\ref{thmtightness-posterior}, that the asymptotic variances of the
bootstrap particle predictor and filter approximations are tight for
any initial distribution $\chi$ such that $\chi(\set{D}) > 0$.

\section{Proofs}
\label{secproofs}
\subsection{Forgetting of the initial distribution}
%
\begin{lemma} \label{lemforgetting-DDelta}
Assume \textup{\hyperref[assum:likelihoodDrift]{(A1)}--\hyperref[assum:majo-g]{(A2)}}. Then for all $\gamma>
2/3$ there exist functions $\rho_\gamma\dvtx  \ooint{0, 1}\ \to\ \ooint{0,1}$
and $C_\gamma\dvtx  \ooint{0,1}\ \to \rset_+$ such that for all $n \in\nset$
and all $\chunk{z}{0}{n-1} \in\set{Y}^{nr}$, where $r \in\nsetpos$ is
as in \hyperref[assum:likelihoodDrift]{\textup{(A1)}} and $z_i = \chunk{y}{ir}{(i+1)r-1}$, satisfying
\[
n^{-1} \sum_{i = 0}^{n - 1}
\1_{\set{K}}(z_i) \geq\gamma,
\]
all functions $f$ and $h$ in $\pmf{\field{X}}$, all finite measures
$\chi$ and $\chi'$ in $\meas{\field{X}}$, and all $\eta\in\,\ooint{0, 1}$,
%
\begin{eqnarray}
&& \bigl\llvert \DDeltaantri{\chi, \chi'} {\chunk{z} {0}
{n-1}} {f,h} \bigr\rrvert\nonumber
\\
&&\qquad \leq\rho^n_\gamma(\eta) \bigl( \chi\Lkhantri{\chunk{z} {0} {n -
1}} f \times\chi' \Lkhantri{\chunk{z} {0} {n - 1}} h +
\chi' \Lkhantri{\chunk{z} {0} {n - 1}} f \times\chi\Lkhantri{\chunk{z} {0} {n
- 1}} h \bigr)\label{itemzfirst}
\\
&&\quad\qquad{}+ C_\gamma(\eta) \eta^n \supnorm{f} \supnorm{h} \prod
_{i = 0}^{n - 1} \supnormantri[2]{
\Lkh{z_i}\1_{\set{X}}} \chi(\set{X}) \chi'(
\set{X}),\nonumber
\\
&& \biggl\llvert \ln \biggl( \frac{\chi\Lkh{\chunk{z}{0}{n-1}}h}{\chi\Lkh{\chunk
{z}{0}{n-1}}f} \biggr) - \ln
\biggl( \frac{\chi' \Lkh{\chunk{z}{0}{n-1}}h}{\chi' \Lkh{\chunk
{z}{0}{n-1}}f} \biggr) \biggr\rrvert\nonumber
\\
\label{itemzsecond}  &&\qquad  \leq\bigl(1-\rho_\gamma(
\eta)\bigr)^{-1}
\\
&&\quad\qquad{}\times \biggl( 2 \rho^n_\gamma(\eta) + \frac{C_\gamma(\eta) \eta^n \supnorm{f} \supnorm{h} \prod_{i=0}^{n-1}
\supnorm{\Lkh{z_i}\1_{\set{X}}} \chi(\set{X})
\chi'(\set{X})}{\chi\Lkh{\chunk{z}{0}{n-1}} f   \times  \chi' \Lkh
{\chunk{z}{0}{n-1}} h}
\biggr),\hspace*{-25pt}\nonumber
\\%
&& \biggl\llvert \frac{\chi\Lkh{\chunk{z}{0}{n-1}}h}{\chi\Lkh{\chunk{z}{0}{n-1}}f} - \frac{\chi' \Lkh{\chunk{z}{0}{n-1}}h}{\chi' \Lkh{\chunk{z}{0}{n-1}}f} \biggr
\rrvert \nonumber
\\
\label{itemzthird}  &&\qquad \leq\rho_\gamma^n(\eta) \biggl( \frac{\chi\Lkh{\chunk
{z}{0}{n-1}}h}{\chi\Lkh{\chunk{z}{0}{n-1}}f} +
\frac{\chi' \Lkh{\chunk
{z}{0}{n-1}}h}{\chi' \Lkh{\chunk{z}{0}{n-1}}f} \biggr)
\\
&&\quad\qquad{} + \frac{C_\gamma(\eta) \eta^n \supnorm{h} \supnorm{f} \prod_{i=0}^{n-1}
\supnorm[2]{\Lkh{z_i} \1_{\set{X}}} \chi(\set{X}) \chi'(\set{X})} {
\chi\Lkh{\chunk{z}{0}{n-1}} f   \times  \chi' \Lkh{\chunk
{z}{0}{n-1}} f}.\nonumber
\end{eqnarray}
\end{lemma}

\begin{pf}
The proof is adapted straightforwardly from \cite{doucmoulines2011}, Proposition~5.
\end{pf}

%
\begin{lemma} \label{lemk+m}
Assume \textup{\hyperref[assum:likelihoodDrift]{(A1)}}. Then there exists a constant $\kappa
> 0$ such that for all $\chi\in\mdr$ [where $\mdr$ is defined in (\ref{eqmort-de-rire})],
%
\begin{equation}
\label{eqmino-loglkh-kplusm} \inf_{(k, m) \in\nsetpos\times\nset} \kappa^{(k + m)} \chi\Lkhantri {
\chunk{Y} {-m} {k-1}} \1_\set{X}> 0,  \qquad\probalone\mbox{-a.s.},
\end{equation}
and
%
\begin{equation}
\label{eqmino-loglkh-kplusplusm} \inf_{(k, m) \in\nsetpos\times\nset} \kappa^{(k + m)} \supnormTxtantri{
\Lkhantri{\chunk{Y} {-m} {k-1}} \1_\set{X}} > 0,  \qquad\probalone\mbox{-a.s.}
\end{equation}
\end{lemma}
\begin{pf}
To derive \eqref{eqmino-loglkh-kplusm} we first establish that
%
\begin{eqnarray}
\label{eqlower-bound-kplusm}
&& \liminf_{k + m \to\infty} (k + m)^{-1} \ln\chi
\Lkhantri{\chunk {Y} {-m} {k - 1}} \1_\set{X}
\nonumber\\[-8pt]\\[-8pt]
&&\qquad \geq-r \espBig{\ln^- \inf_{x \in\set{D}} \delta_x \Lkhantri{
\chunk{Y} {0} {r - 1}} \1_\set{D}} > - \infty,  \qquad \probalone\mbox{-a.s.},\nonumber
\end{eqnarray}
where the last inequality follows from \hyperref[assum:likelihoodDrift]{(A1)}(iii). We now establish the first inequality in
\eqref{eqlower-bound-kplusm}. Set $a_{k, m} \eqdef-k + \lfloor
(k+m)/r \rfloor r$ and note that $- a_{k, m} \in\{-m, \ldots, - m + r
- 1 \}$. Then write
%
\begin{eqnarray}\label{eqonestep}
&&\ln\chi\Lkhantri{\chunk{Y} {-m} {k-1}} \1_\set{X} \nonumber
\\
&&\qquad \geq \ln\chi\Lkhantri{\chunk{Y} {-m} {-a_{k,m}}}
\1_\set{D} + \sum_{i = 0}^{\lfloor(k + m)/r \rfloor-1} \ln
\inf_{x \in\set{D}} \delta _x \Lkhantri{\chunk{Y} {-
a_{k, m} + i r} {-a_{k, m} + (i + 1)r -1}} \1_\set
{D}
\nonumber\\[-8pt]\\[-8pt]
&&\qquad \geq- \sum_{i = 0}^{r - 1} \ln^- \chi
\Lkhantri{\chunk{Y} {-m} {-m + i}} \1_\set{D}\nonumber
\\
&&\quad\qquad{}  - \sum
_{i = 0}^{\lfloor(k+m)/r \rfloor-1} \ln^- \inf_{x
\in\set{D}}\delta_x\Lkhantri{\chunk{Y} {-a_{k, m} + i r} {-a_{k, m}
+ (i + 1)r -1}} \1_\set{D}.\nonumber
\end{eqnarray}
For $i \in\nset$, set $\modulotxt{i} \eqdef i - \lfloor i / r \rfloor
r$. With this notation, $a_{k, m}= \modulotxt{a_{k,m}}+\lfloor a_{k,m}
/ r \rfloor r $.
Then, since $\modulotxt{i} \in\{0,\ldots,r-1\}$,
%
\begin{eqnarray}
\label{eqtwostep} && - \sum_{i = 0}^{\lfloor(k + m)/r \rfloor- 1}
\ln^- \inf_{x
\in\set{D}} \delta_x \Lkhantri{\chunk{Y}
{-a_{k, m} + i r} {-a_{k, m} + (i + 1)r -1}}
\1_\set{D}\nonumber
\\
&&\qquad = - \sum_{i = 0}^{\lfloor(k + m)/r \rfloor-1} \ln^- \inf
_{x \in\set
{D}} \delta_x \Lkhantri{\chunk{Y} {-
\modulotxt{a_{k,m}}+\bigl(i - \lfloor a_{k, m} / r \rfloor
\bigr) r} {-\modulotxt{a_{k,m}} + \bigl(i - \lfloor a_{k, m} /
r \rfloor+ 1\bigr)r -1 }} \1_\set{D}
\nonumber\\[-8pt]\\[-8pt]
&&\qquad \geq- \sum_{j = 0}^{r - 1} \sum
_{i = 0}^{\lfloor(k + m) / r \rfloor
- 1} \ln^- \inf_{x \in\set{D}}
\delta_x \Lkhantri{\chunk{Y} {-j + \bigl(i - \lfloor a_{k, m} / r
\rfloor\bigr) r} {-j + \bigl(i - \lfloor a_{k,m} / r \rfloor +1\bigr)r
-1}} \1_\set{D}
\nonumber
\\
&&\qquad = - \sum_{j = 0}^{r - 1} \sum
_{\ell=- \lfloor a_{k,m} / r \rfloor
}^{\lfloor(k+m) / r \rfloor- \lfloor a_{k,m} / r \rfloor- 1} \ln^- \inf_{x \in\set{D}}
\delta_x\Lkhantri{\chunk{Y} {-j+\ell r} {-j+(\ell+1)r -1}}
\1_\set{D},\nonumber
\end{eqnarray}
where the last identity follows by reindexing the summation.
We now plug \eqref{eqtwostep} into \eqref{eqonestep}; the ergodicity
of the process $\sequence{Z}$ [assumption \hyperref[assum:likelihoodDrift]{(A1)}(i)] then implies, via
Lemma~\ref{lemstationarykplusm}, \probalone\mbox{-a.s.},
\begin{eqnarray*}
&& \liminf_{k+m \to\infty} (k+m)^{-1} \ln\chi\Lkhantri{\chunk {Y}
{-m} {k-1}} \1_\set{X}
\\
&&\qquad \geq\sum_{j = 0}^{r - 1} \espBig{\ln^- \inf
_{x \in\set{D}} \delta_x \Lkhantri {\chunk{Y} {-j} {-j + r - 1}}
\1_\set{D}} = - r \espBig{\ln^- \inf_{x \in
\set{D}}
\delta_x \Lkhantri{\chunk{Y} {0} {r-1}} \1_\set{D}},
\end{eqnarray*}
which shows \eqref{eqlower-bound-kplusm}. Now, choose a constant
$\kappa$ such that
\[
-r \espBig{\ln^- \inf_{x \in\set{D}} \delta_x \Lkhantri{\chunk{Y}
{0} {r-1}} \1 _\set{D}} > - \ln\kappa > -\infty.
\]
According to \eqref{eqlower-bound-kplusm}, there exists a \probalone
\mbox{-a.s.} finite $\nsetpos$-valued random variable $N$ such that if $k + m
\geq N$,
\[
\ln\chi\Lkhantri{\chunk{Y} {0} {r - 1}} \1_{\set{X}} \geq(- \ln\kappa) (k + m),
\]
which implies that
\[
\inf_{k+m \geq N} \kappa^{k+m} \chi\Lkhantri{\chunk{Y} {0}
{r-1}} \1_{\set
{X}} \geq1.
\]
On the other hand, assumption~\hyperref[assum:majo-g]{(A2)} implies that for all
$(k, m) \in\nsetpos\times\nset$, $\chi\Lkh{\chunk{Y}{0}{r - 1}} \1
_{\set{X}} > 0$, $\probalone\mbox{-a.s.}$ This completes the proof of \eqref
{eqmino-loglkh-kplusm}. Finally, the proof of~\eqref{eqmino-loglkh-kplusplusm} follows by combining
\[
\supnormTxtantri{\Lkhantri{\chunk{Y} {-m} {k-1}} \1_\set{X}} \geq\chi\Lkhantri{
\chunk {Y} {-m} {k-1}} \1_\set{X} %
\]
and \eqref{eqmino-loglkh-kplusm}.
\end{pf}

For all probability measures $\chi\in\probmeas{\field{X}}$, all $(k,
m) \in\nsetpos\times\nset$, and all sequences $\chunk{y}{-m}{k} \in
\set{Y}^{m+k+1}$, define the set
%
\begin{eqnarray}
\label{eqdefinition-Pkm}
&& \mcpantri{\chunk{y} {-m} {k}} {\chi}
\nonumber\\[-8pt]\\[-8pt]
&&\qquad \eqdef \bigl\{ \tilde{\chi} \in\probmeas{\field{X}} \dvtx  \supnormantri{\llhd
{y_k}} \times\tilde{\chi} \Lkhantri{\chunk{y} {-m} {k-1}}
\1_{\set{X}} \geq (1/2) \chi\Lkhantri{\chunk{y} {-m} {k}} \1_{\set{X}}
\bigr\}\nonumber
\end{eqnarray}
of probability measures on $(\set{X}, \field{X})$
and note that this set is nonempty since $\chi\in\mcp{\chunk
{y}{-m}{k}}{\chi}$. The choice of $1/2$ in the definition of $\mcp{\chunk
{y}{-m}{k}}{\chi}$ is irrelevant, and this factor can be replaced by
any constant strictly less than~$1$.

\begin{proposition} \label{propbound-difference-improved}
Assume \textup{\hyperref[assum:likelihoodDrift]{(A1)}--\hyperref[assum:majo-g]{(A2)}}. Then there exists a
constant $\beta\in\,\ooint{0,1}$ such that the following hold:
\begin{longlist}[(iii)]
\item[(i)]\label{itemupperbound-log-lkh-diff-init-distr}
For all probability measures $\chi$ and $\chi'$ in $\mdr$ there exists
a $\probalone$\mbox{-a.s.} finite random variable $C_{\chi,\chi'}$ such that
for all $(k, m) \in\nsetpos\times\nset$ and all $\tilde{\chi} \in
\mcp{\chunk{Y}{-m}{k}}{\chi}$,
\[
\ln \biggl( \frac{\tilde{\chi} \Lkh{\chunk{Y}{-m}{k}} \1_\set{X}}{\tilde
{\chi} \Lkh{\chunk{Y}{-m}{k-1}}\1_\set{X}} \biggr) - \ln \biggl( \frac{\chi' \Lkh{\chunk{Y}{-m}{k}} \1_\set{X}}{\chi' \Lkh{\chunk
{Y}{-m}{k-1}}\1_\set{X}} \biggr)
\leq C_{\chi,\chi'} \beta^{k+m},  \qquad\probalone\mbox{-a.s.}
\]
\item[(ii)]\label{itemupperbound-log-lkh-one-step}
For all probability measures $\chi$ in $\mdr$ there exists a
$\probalone$\mbox{-a.s.} finite random variable $C_{\chi}$ such that for all
$(k, m) \in\nsetpos\times\nset$,
\[
\biggl\llvert \ln \biggl( \frac{\chi\Lkh{\chunk{Y}{-m}{k}} \1_\set{X}}{\chi
\Lkh{\chunk{Y}{-m}{k-1}}\1_\set{X}} \biggr) - \ln \biggl(
\frac{\chi\Lkh{\chunk{Y}{-m-1}{k}} \1_\set{X}}{\chi\Lkh
{\chunk{Y}{-m-1}{k-1}}\1_\set{X}} \biggr) \biggr\rrvert \leq C_\chi
\beta^{k + m},  \qquad\probalone\mbox{-a.s.}
\]
\item[(iii)]\label{itemevry2}
There exists a $\probalone$\mbox{-a.s.} finite random variable $C$ such that
for $m \in\nsetpos$, all probability measures $\chi$ and $\chi'$ in
$\probmeas{\field{X}}$ and all $h \in\bmf{\field{X}}$,
\[
\frac{|\DDelta{\chi,\chi'}{\chunk{Y}{-m}{-1}}{h, \1_{\set{X}}}|}{
\| \Lkh{\chunk{Y}{-m}{-1}} \1_{\set{X}} \|_{\infty}^2} \leq C \beta^m \supnorm{h},  \qquad\probalone\mbox{-a.s.}
\]
\end{longlist}
\end{proposition}

\begin{pf*}{Proof of Proposition~\ref{propbound-difference-improved}\hyperref[itemupperbound-log-lkh-diff-init-distr]{\textup{(i)}}
and \hyperref[itemupperbound-log-lkh-one-step]{\textup{(ii)}}}
Let $\tilde{\chi} \in\mcp{\chunk{Y}{-m}{k}}{\chi}$.
Recall the notation $Z_i = \chunk{Y}{ir}{(i + 1)r - 1}$ and consider
the decompositions
\begin{eqnarray*}
\chi\Lkhantri{\chunk{Y} {-m} {k}} \1_\set{X} &=& \chi\Lkhantri{\chunk{Y} {-m} {-
\lfloor m / r \rfloor r - 1}} \Lkhantri{\chunk{Z} {- \lfloor m / r \rfloor } {\lfloor
k / r \rfloor-1}} \Lkhantri{\chunk{Y} {\lfloor k / r \rfloor r} {k}}
\1_\set{X},
\\
\chi\Lkhantri{\chunk{Y} {-m} {k-1}} \1_\set{X} &=& \chi\Lkhantri{\chunk {Y} {-m}
{-\lfloor m / r \rfloor r - 1}} \Lkhantri{\chunk{Z} {-\lfloor m / r \rfloor} {\lfloor
k / r \rfloor-1}} \Lkhantri{\chunk{Y} {\lfloor k / r \rfloor r} {k - 1}}
\1_\set{X},
\end{eqnarray*}
where we make use of convention \eqref{eqconvention} if necessary.

Choose $\gamma$ such that $2/3 < \gamma< \probTxt{Z_0 \in\set{K}}$,
where $\set{K}$ is defined in \hyperref[assum:likelihoodDrift]{(A1)}(i). Assume that $(k, m) \in\nsetpos\times\nset$
are both larger than $r$ and denote by $b_{k, m} \eqdef\lfloor k / r
\rfloor+ \lfloor m / r \rfloor$. In addition, define the event
\[
\label{eqironarm} \Omega_{k, m} \eqdef \Biggl\{ \biggl( \biggl\lfloor
\frac{k}{r} \biggr\rfloor+ \biggl\lfloor\frac{m}{r} \biggr\rfloor
\biggr)^{-1} \sum_{\ell= - \lfloor m / r \rfloor}^{ \lfloor k / r \rfloor- 1}
\1_{\set
{K}}(Z_\ell) \geq\gamma \Biggr\}.
\]
By Lemma~\ref{lemforgetting-DDelta} [equation~\eqref
{itemzsecond}] it holds for all $\eta\in\,\ooint{0,1}$, on the event
$\Omega_{k, m}$,
%
\begin{eqnarray}
\label{eqjuvisy} && \bigl(1-\rho_\gamma(\eta)\bigr) \biggl( \ln
\biggl( \frac{\tilde{\chi} \Lkh
{\chunk{Y}{-m}{k}} \1_\set{X}}{\tilde{\chi} \Lkh{\chunk{Y}{-m}{k - 1}}
\1_\set{X}} \biggr)- \ln \biggl( \frac{\chi' \Lkh{\chunk{Y}{-m}{k}} \1
_\set{X}}{\chi' \Lkh{\chunk{Y}{-m}{k-1}}\1_\set{X}} \biggr)
\biggr)\nonumber
\\
&&\qquad \stackrel{(\mathrm{a})} {\leq} 2 \rho^{b_{k,m}}_\gamma(
\eta) + \frac{C_\gamma(\eta) \eta^{b_{k,m}} \supnorm{\llhd{Y_k}} \prod_{i=-m}^{k-1} \supnorm[2]{\llhd{Y_i}}}{\tilde{\chi} \Lkh{\chunk
{Y}{-m}{k-1}}\1_\set{X}   \times  \chi' \Lkh{\chunk{Y}{-m}{k}}\1_\set
{X}}
\\
&&\qquad\stackrel{(\mathrm{b})} {\leq} 2 \rho^{b_{k,m}}_\gamma(\eta) +
\frac{2 C_\gamma(\eta) \eta^{b_{k, m}} \prod_{i=-m}^k \supnorm
[2]{\llhd{Y_i}}}{\chi\Lkh{\chunk{Y}{-m}{k}}\1_\set{X}   \times  \chi
' \Lkh{\chunk{Y}{-m}{k}}\1_\set{X}},
\nonumber
\end{eqnarray}
where:
\begin{longlist}[(a)]
\item[(a)] follows from \eqref{itemzsecond} and the bound $\delta_x \Lkh
{\chunk{Y}{u}{v}} \1_\set{X} \leq\prod_{\ell=u}^v \supnorm{\llhd{Y_\ell
}}$, valid for $u \leq v$, and
\item[(b)] follows from the fact that $\tilde{\chi} \in\mcp{\chunk
{Y}{-m}{k}}{\chi}$.
\end{longlist}
Since, under \hyperref[assum:likelihoodDrift]{(A1)}(i), the
sequence $\sequence{Z}$ is ergodic and $\probTxt{Z_0 \in\set{K}} >
\gamma$,
Lemma~\ref{lemstationarykplusm} implies that
\[
\probBig{\bigcup_{j \geq0} \
\mathop{\bigcap_{(k,m) \in
\nsetpos\times\nset}}_{k + m \geq j}
\Omega_{k, m}}=1.  %
\]
Hence, there exists a \probalone\mbox{-a.s.} finite integer-valued random
variable $U$ such that \eqref{eqjuvisy} is satisfied for all $(k, m)
\in\nsetpos\times\nset$ such that $k + m \geq U$.

The lower bound obtained in Lemma~\ref{lemk+m} implies that there
exists a constant $\kappa> 0$ such that for all probability measures
$\chi$ and $\chi'$ in $\mdr$ and all $(k,m) \in\nsetpos\times\nset$,
\probalone\mbox{-a.s.},
\begin{eqnarray*}
\chi\Lkhantri{\chunk{Y} {-m} {k}} \1_\set{X} &\geq&\widebar{C}_{\chi,\chi'}
\kappa ^{-(k+m+1)},
\\
\chi' \Lkhantri{\chunk{Y} {-m} {k}} \1_\set{X} &\geq&
\widebar{C}_{\chi,\chi'} \kappa^{-(k+m+1)},
\end{eqnarray*}
where $\widebar{C}_{\chi,\chi'}$ is a \probalone\mbox{-a.s.} finite constant.

By plugging these bounds into \eqref{eqjuvisy} and using
Lemma~\ref{lemboundFilter} with
$\eta$ sufficiently small (note that \eqref{eqjuvisy} is satisfied for
all $\eta\in\,\ooint{0,1}$),
we conclude that there exist a $\probalone\mbox{-a.s.}$ finite random variable
$C_{\chi, \chi'}$
and a constant $\beta< 1$ such that for all $(k,m) \in\nsetpos\times
\nset$, \probalone\mbox{-a.s.},
\[
\ln \biggl( \frac{\tilde{\chi} \Lkh{\chunk{Y}{-m}{k}} \1_\set{X}}{\tilde
{\chi} \Lkh{\chunk{Y}{-m}{k-1}}\1_\set{X}} \biggr)- \ln \biggl( \frac{\chi
' \Lkh{\chunk{Y}{-m}{k}} \1_\set{X}}{\chi' \Lkh{\chunk{Y}{-m}{k-1}}\1
_\set{X}} \biggr)
\leq C_{\chi,\chi'} \beta^{k+m},  %
\]
which completes the proof of Proposition~\ref{propbound-difference-improved}\hyperref[itemupperbound-log-lkh-diff-init-distr]{(i)}. Note that $\chi\in\mcp
{\chunk{Y}{-m}{k}}{\chi}$ implies that the previous relation is
satisfied with $\tilde{\chi}= \chi$.

The proof of Proposition \ref{propbound-difference-improved}\hyperref[itemupperbound-log-lkh-one-step]{\textup{(ii)}} follows the same
lines as the proof of Proposition~\ref{propbound-difference-improved}\hyperref[itemupperbound-log-lkh-diff-init-distr]{(i)}
and is omitted for brevity.
\end{pf*}

\begin{pf*}{Proof of Proposition \ref{propbound-difference-improved}\hyperref[itemevry2]{\textup{(iii)}}}
We may assume that the function $h$ is nonnegative (otherwise the positive and negative parts of $h$ can be treated separately).
As in the proof of Proposition~\ref{propbound-difference-improved}\hyperref[itemupperbound-log-lkh-diff-init-distr]{(i)}, write
\[
\chi\Lkhantri{\chunk{Y} {-m} {-1}} h = \chi\Lkhantri{\chunk{Y} {-m} {-\lfloor m / r
\rfloor r-1}} \Lkhantri{\chunk{Z} {-\lfloor m / r \rfloor} {-1}} h %
\]
and define the event
\[
\label{eqironarmtwo} \Omega_m \eqdef \Biggl\{ \biggl\lfloor
\frac{m}{r} \biggr\rfloor^{-1} \sum_{\ell=- \lfloor m / r  \rfloor}^{ -1}
\1_{\set{K}}(Z_\ell ) \geq\gamma \Biggr\}.
\]
By Lemma~\ref{lemforgetting-DDelta} [equation~(\ref{itemzthird})]
it holds, on the event $\Omega_m$,
%
\begin{eqnarray}
\label{eqvertdemaison}
&& \biggl\llvert \frac{\chi\Lkh{\chunk{Y}{-m}{-1}} h}{\chi\Lkh{\chunk
{Y}{-m}{-1}}\1_\set{X}} - \frac{\chi' \Lkh{\chunk{Y}{-m}{-1}} h}{\chi' \Lkh{\chunk
{Y}{-m}{-1}}\1_\set{X}} \biggr\rrvert
\nonumber\\[-8pt]\\[-8pt]
&&\qquad \leq2 \supnorm{h} \rho_\gamma^{\lfloor m / r \rfloor}(\eta) + \frac
{C_\gamma(\eta) \eta^{\lfloor m / r \rfloor} \supnorm{h} \prod_{i=-m}^{-1} \supnorm[2]{\llhd{Y_i}}} {
\chi\Lkh{\chunk{Y}{-m}{-1}}\1_\set{X}   \times  \chi' \Lkh{\chunk
{Y}{-m}{-1}}\1_\set{X}},\nonumber
\end{eqnarray}
where we used that for $u\leq v$, $\delta_x \Lkh{\chunk{Y}{u}{v}} \1
_\set{X} \leq\prod_{\ell= u}^v \supnorm{\llhd{Y_\ell}}$. Under \hyperref[assum:likelihoodDrift]{(A1)}(i), Birkhoff's ergodic
theorem (see, e.g., \cite{shiryaev1996}) ensures that
$\probTxt{\liminf_{m \to\infty} \Omega_m }= 1$; therefore, there
exists a $\probalone\mbox{-a.s.}$ finite random variable $U$ such that \eqref
{eqvertdemaison} is satisfied for $m \geq U$. Then, for $m \geq U$,
%
\begin{eqnarray}
\label{eqrisorangis}
&& \frac{|\DDelta{\chi,\chi'}{\chunk{Y}{-m}{-1}}{h, \1_{\set{X}}}|}{
\| \Lkh{\chunk{Y}{-m}{-1}} \1_{\set{X}}\|^2}\nonumber
\\
&&\qquad =\frac{\chi\Lkh{\chunk{Y}{-m}{-1}}\1_\set{X}   \times  \chi
' \Lkh{\chunk{Y}{-m}{-1}}\1_\set{X}}{\| \Lkh{\chunk{Y}{-m}{-1}} \1_{\set
{X}} \|^2 } \biggl\llvert \frac{\chi\Lkh{\chunk{Y}{-m}{-1}} h}{\chi\Lkh
{\chunk{Y}{-m}{-1}}\1_\set{X}} -
\frac{\chi' \Lkh{\chunk{Y}{-m}{-1}}
h}{\chi' \Lkh{\chunk{Y}{-m}{-1}}\1_\set{X}} \biggr\rrvert
\\
&&\qquad \leq2 \supnorm{h} \rho_\gamma^{\lfloor m / r \rfloor}(\eta) +
\frac{C_\gamma(\eta) \eta^{\lfloor m / r \rfloor} \supnorm{h} \prod_{i=-m}^{-1} \supnorm[2]{\llhd{Y_i}}}{\| \Lkh{\chunk{Y}{-m}{-1}} \1
_{\set{X}} \|^2},\nonumber
\end{eqnarray}
we have used that $\chi\Lkh{\chunk{Y}{-m}{-1}} \1_\set{X} \leq
\supnormTxt{\Lkh{\chunk{Y}{-m}{-1}} \1_{\set{X}}}$. By Lemma~\ref{lemk+m}, equation~\eqref{eqmino-loglkh-kplusplusm}, there exist a
constant $\kappa> 0$ and a $\probalone$\mbox{-a.s.} finite random variable
$\widebar{C}$ such that
\[
\supnormTxtantri{\Lkhantri{\chunk{Y} {-m} {-1}} \1_{\set{X}}} \geq\widebar{C} \kappa
^{-m},  \qquad \probalone\mbox{-a.s.} %
\]
Finally, we complete the proof by inserting this bound into \eqref
{eqrisorangis} and applying Lemma~\ref{lemboundFilter} to the
right-hand side of the resulting inequality.
\end{pf*}
%

\subsection{Convergence of the log-likelihood}
%
\begin{lemma} \label{lemmekong}
Assume \textup{\hyperref[assum:likelihoodDrift]{(A1)}--\hyperref[assum:majo-g]{(A2)}}. Then, $\probalone$\mbox{-a.s.},
%
\begin{eqnarray}
\label{eqcvgceforward} \lim_{n \to\infty} n^{-1} \ln\supnormantri{
\Lkhantri{\chunk{Y} {0} {n}} \1_\set {X}} &=& \ell_\infty,
\\
\label{eqcvgcebackward}  \lim_{n \to\infty} n^{-1} \ln\supnormTxtantri{
\Lkhantri{\chunk{Y} {-n} {0}} \1 _\set{X}} &=& \ell_\infty,
\\
\label{eqcvgcebackward-stat}  \lim_{n \to\infty} n^{-1} \sum
_{k=1}^n \ln\likeliantri{\chunk{Y} {-\infty } {-k-1}}
{Y_{-k}}&=& \ell_\infty,
\end{eqnarray}
where $\ell_\infty$ is defined in \eqref{eqdefinitionellinfty}.
\end{lemma}

\begin{pf*}{Proof of \eqref{eqcvgceforward}}
Let $(\alpha_n)_{n \in\nsetpos}$ be a nondecreasing sequence such
that $\lim_{n \to\infty} \alpha_n = 1$ and for any $n \in\nsetpos$,
$\alpha_n \geq1/2$.
For all $n \in\nset$, choose $\tilde{x}_n \in\set{X}$ such that
%
\begin{equation}
\label{eqdefinition-alpha} \alpha_n \supnormantri{\Lkhantri{\chunk{Y} {0} {n}}
\1_\set{X}} \leq\delta_{\tilde
{x}_n} \Lkhantri{\chunk{Y} {0} {n}}
\1_\set{X} \leq\supnormantri{\Lkhantri{\chunk {Y} {0} {n}} \1_\set{X}}.
\end{equation}
Note that for all $k \in\nsetpos$,
%
\begin{equation}
\label{eqrel1} \delta_{\tilde{x}_{k-1}} \Lkhantri{\chunk{Y} {0} {k-1}}
\1_\set{X} \geq\alpha _{k-1} \supnormTxtantri{ \Lkhantri{\chunk{Y} {0}
{k-1}} \1_\set{X}} \geq\alpha _{k-1} \delta_{\tilde{x}_k}
\Lkhantri{\chunk{Y} {0} {k-1}} \1_\set{X}.
\end{equation}
On the other hand, for all probability measures $\chi\in\probmeas{\field{X}}$ it holds that
%
\begin{equation}
\label{eqrel2} \delta_{\tilde{x}_k} \Lkhantri{\chunk{Y} {0} {k-1}}
\1_\set{X} \stackrel {\mathrm{(a)}} {\geq} \frac{\delta_{\tilde{x}_k} \Lkh{\chunk{Y}{0}{k}} \1_\set
{X}}{\supnormTxt{\llhd{Y_k}}} \stackrel{\mathrm{(b)}} {
\geq} \alpha_k \frac{\supnormTxt{\Lkh{\chunk{Y}{0}{k}}
\1_\set{X}}}{\supnormTxt{\llhd{Y_k}}} \geq\alpha_k
\frac{\chi\Lkh{\chunk{Y}{0}{k}}\1_\set{X}}{\supnormTxt
{\llhd{Y_k}}},
\end{equation}
where (a) follows from the bound $\delta_{\tilde{x}_k} \Lkh{\chunk
{Y}{0}{k}} \1_{\set{X}} \leq\supnorm{\llhd{Y_k}} \delta_{\tilde{x}_k}
\Lkh{\chunk{Y}{0}{k-1}} \1_{\set{X}}$ and (b) stems from the definition
\eqref{eqdefinition-alpha} of $\alpha_n$.
Then
%
\begin{eqnarray}\label{eqbound-difference-4}
0 &\leq& n^{-1} \bigl( \ln \supnormantri{ \Lkhantri{\chunk{Y} {0} {n}}
\1_\set{X}} - \ln\chi\Lkhantri{\chunk{Y} {0} {n}} \1_\set{X} \bigr)
\nonumber
\\
& \leq&-n^{-1} \ln\alpha_n + n^{-1} \bigl( \ln \bigl(
\alpha_n \supnormantri{ \Lkhantri{\chunk{Y} {0} {n}} \1_\set{X}} \bigr) -
\ln\chi\Lkhantri {\chunk{Y} {0} {n}} \1_\set{X} \bigr)
\nonumber
\\
& \leq&-n^{-1} \ln\alpha_n + n^{-1} \bigl( \ln
\delta_{\tilde{x}_n} \Lkhantri{\chunk{Y} {0} {n}} \1_\set{X} - \ln\chi
\Lkhantri{\chunk{Y} {0} {n}} \1 _\set{X} \bigr)
\\
& =& -n^{-1} \ln\alpha_n + n^{-1} \bigl( \ln
\delta_{\tilde{x}_0} \Lkh {Y_0} \1_\set{X} - \ln\chi
\Lkh{Y_0} \1_\set{X} \bigr)
\nonumber
\\
&&{} + n^{-1} \sum_{k=1}^n
\biggl[ \ln \biggl( \frac{\delta_{\tilde{x}_k} \Lkh{\chunk{Y}{0}{k}} \1_\set{X}}{\delta
_{\tilde{x}_{k-1}} \Lkh{\chunk{Y}{0}{k-1}} \1_\set{X}} \biggr) - \ln \biggl( \frac{\chi\Lkh{\chunk{Y}{0}{k}}\1_\set{X}}{\chi\Lkh
{\chunk{Y}{0}{k-1}} \1_\set{X}}
\biggr) \biggr].\nonumber
\end{eqnarray}
For each term in the sum it holds, by \eqref{eqrel1},
\begin{eqnarray*}
&& \ln \biggl( \frac{\delta_{\tilde{x}_k} \Lkh {\chunk{Y}{0}{k}} \1_\set{X}}{\delta
_{\tilde{x}_{k-1}} \Lkh {\chunk{Y}{0}{k-1}} \1_\set{X}} \biggr) - \ln \biggl( \frac{\chi\Lkh {\chunk{Y}{0}{k}}\1_\set{X}}{\chi\Lkh
{\chunk{Y}{0}{k-1}} \1_\set{X}}
\biggr)
\nonumber
\\
&&\qquad \leq- \ln\alpha_{k-1} + \ln \biggl( \frac{\delta_{\tilde{x}_k} \Lkh
{\chunk{Y}{0}{k}}\1_\set{X}}{\delta_{\tilde{x}_k} \Lkh {\chunk
{Y}{0}{k-1}} \1_\set{X}} \biggr) - \ln
\biggl( \frac{\chi\Lkh {\chunk{Y}{0}{k}}\1_\set{X}}{\chi\Lkh {\chunk
{Y}{0}{k-1}}\1_\set{X}} \biggr).
\end{eqnarray*}
For all $k \in\nsetpos$, \eqref{eqrel2} implies that
\[
\delta_{\tilde{x}_k} \Lkhantri{\chunk{Y} {0} {k-1}} \1_\set{X} \geq
\frac{1}2 \frac{\chi\Lkh {\chunk{Y}{0}{k}}\1_\set{X}}{\supnorm{\llhd{Y_k}}},
\]
so that $\delta_{\tilde{x}_k}$ belongs to the set $\mcp{\chunk
{Y}{0}{k-1}}{\chi}$ [defined in~\eqref{eqdefinition-Pkm}].
Proposition~\ref{propbound-difference-improved}\hyperref[itemupperbound-log-lkh-diff-init-distr]{(i)} then provides a constant
$\beta\in\,\ooint{0,1}$ and a \probalone\mbox{-a.s.} finite
random variable $C_{\chi}$ such that
%
\begin{equation}
\ln \biggl( \frac{\delta_{\tilde{x}_k} \Lkh {\chunk{Y}{0}{k}}\1_\set
{X}}{\delta_{\tilde{x}_k} \Lkh {\chunk{Y}{0}{k-1}} \1_\set{X}} \biggr) - \ln \biggl( \frac{\chi\Lkh {\chunk{Y}{0}{k}}\1_\set{X}}{\chi\Lkh {\chunk
{Y}{0}{k-1}}\1_\set{X}} \biggr)
\leq C_\chi\beta^k.  \label
{eqbound-difference-5}
\end{equation}
Finally, statement \eqref{eqcvgceforward} follows by plugging the
bound \eqref{eqbound-difference-5} into \eqref{eqbound-difference-4},
letting $n$ tend to infinity and using \eqref{itemlim-ponct-rate-entropy-Y}.
\end{pf*}

\begin{pf*}{Proof of \eqref{eqcvgcebackward}}
For all $(p, n) \in\nset^2$ such that $p \leq n$, define $W_{p, n}
\eqdef\ln\supnormTxt{\Lkh{\chunk{Y}{p}{n-1}} \1_\set{X}}$ and
$\widetilde{W}_{p, n} \eqdef\ln\supnormTxt{\Lkh{\chunk{Y}{-n+1}{-p}} \1
_\set{X}}$. Note that these two sequences are subadditive in the sense
that for all $(p, n) \in\nset^2$ such that $p \leq n$,
\begin{eqnarray*}
W_{0, n} &\leq& W_{0, p} + W_{p, n},
\\
\widetilde{W}_{0, n} &\leq&\widetilde{W}_{0, p} +
\widetilde{W}_{p, n}.
\end{eqnarray*}
Finally, for all $x \in\set{D}$, $m \in\nset$ and $\chunk{y}{0}{mr -
1} \in\set{Y}^{mr}$, it holds that
%
\begin{equation}
\label{eqangkhor} \supnormTxtantri{\Lkhantri{\chunk{y} {0} {mr-1}} \1_{\set{X}}} \geq
\delta_x \Lkhantri{\chunk{y} {0} {mr - 1}} \1_\set{X} \geq\prod
_{\ell= 0}^{m - 1} \inf_{x \in\set{D}}
\delta_x \Lkhantri{\chunk {y} {kr} {(k+1)r - 1}} \1_\set{D}.
\end{equation}
Using the stationarity of the observation process $\sequence{Y}[k][\zset
]$, we get, via assumption~\hyperref[assum:likelihoodDrift]{(A1)}(iii), for all $m \in\nsetpos$,
%
\begin{eqnarray}
\label{eqlebrasdef}
&& (mr)^{-1} \esp{W_{0, mr}} \nonumber
\\
&&\qquad = (mr)^{-1} \esp{\widetilde{W}_{0, mr}} \geq(mr)^{-1}
\espbig{\ln\supnormTxtantri{\Lkhantri{\chunk{y} {0} {mr-1}} \1_{\set
{X}}} }
\\
&&\qquad \geq r^{-1} \espbig{ \ln\inf_{x \in\set{D}}
\delta_x \Lkhantri{\chunk {y} {kr} {(k+1)r - 1}} \1_\set{D} } >
- \infty.\nonumber
\end{eqnarray}
The sequences $(\esp{W_{0,n}})_{n \in\nsetpos}$ and
$(\espTxt{\widetilde{W}_{0,n}})_{n \in\nsetpos}$ are subadditive; Fekete's
lemma~(see \cite{polya1976a}) thus implies that
the sequences $(n^{-1} \esp{W_{0,n}})_{n \in\nsetpos}$ and
$(n^{-1} \espTxt{\widetilde{W}_{0,n}})_{n \in\nsetpos}$ have limits in
$\coint{-\infty,\infty}$ and that
\begin{eqnarray*}
\lim_{n \to\infty} n^{-1} \esp{W_{0,n}} &=& \lim
_{n \to\infty} n^{-1} \esp{\widetilde{W}_{0,n}}
\\
&=& \inf_{n \in\nsetpos} n^{-1} \esp{W_{0,n}}
\\
&=& \inf
_{n \in\nsetpos} n^{-1} \esp{\widetilde{W}_{0,n}}.
\end{eqnarray*}
However, by \eqref{eqlebrasdef} there exists a subsequence that is
bounded away from $-\infty$, showing that
\begin{eqnarray*}
\inf_{n \in\nsetpos} n^{-1} \esp{W_{0,n}}&=&\lim
_{n \to\infty} n^{-1} \esp{W_{0,n}} > -\infty,
\\
\inf_{n \in\nsetpos} n^{-1} \esp{\widetilde{W}_{0,n}}&=&
\lim_{n \to\infty
} n^{-1} \esp{\widetilde{W}_{0,n}}
> -\infty.
\end{eqnarray*}
Now,\vspace*{2pt} by applying Kingman's subadditive ergodic theorem (see \cite
{kingman1973}) and using again that $\espTxt{\widetilde{W}_{0, k}} = \espTxt
{W_{0, k}}$ under stationarity, we obtain
\begin{eqnarray*}
\lim_{n \to\infty} n^{-1} \widetilde{W}_{0, n} &=&
\inf_{n \in\nsetpos} n^{-1} \esp{\widetilde{W}_{0, n}}
= \inf_{n \in\nsetpos} n^{-1} \esp{W_{0,
n}}
\\
&=& \lim_{n \to\infty} n^{-1} W_{0, n} =
\ell_\infty,  \qquad\probalone \mbox{-a.s.},
\end{eqnarray*}
where the last limit follows from \eqref{eqcvgceforward}. This
completes the proof of statement~\eqref{eqcvgcebackward}.
\end{pf*}

\begin{pf*}{Proof of \eqref{eqcvgcebackward-stat}}
Since $\espTxt{| \ln\likeli{\chunk{Y}{-\infty}{-1}}{Y_0} |} < \infty$
and the process $\sequence{Y}[k][\zset]$ is stationary and ergodic,
\eqref{eqcvgcebackward-stat} follows from Birkhoff's ergodic theorem.
\end{pf*}
%

\begin{appendix}
\section*{Appendix: Technical lemmas} \label{sectechlemmas}

%
\begin{lemma}\label{lemstationarykplusm}
If $\sequence{U}[n][\zset]$ is a stationary and ergodic sequence of
random variables such that $\esp{|U_0|}<\infty$, then
%
\begin{equation}
\label{eqstationarykplusm} \lim_{k+m \to\infty} (k+m)^{-1} \Biggl(\sum
_{\ell=-m}^{k-1}U_{\ell
} \Biggr) =
\esp{U_0},  \qquad\probalone\mbox{-a.s.}
\end{equation}
\end{lemma}

\begin{pf}
Denote
\begin{eqnarray*}
\Omega_1 &\eqdef& \Biggl\{ \omega\in\Omega; \lim
_{k+m \to\infty} (k+m)^{-1} \Biggl(\sum
_{\ell=-m}^{k-1}U_{\ell}(\omega) \Biggr) =
\esp{U_0} \Biggr\},
\\
\Omega_2 &\eqdef& \biggl\{ \omega\in\Omega; \lim_{m \to\infty}
\frac{\sum_{\ell=-m}^{-1} U_\ell(\omega)}{m}=\lim_{k \to\infty} \frac{\sum_{\ell=0}^{k-1} U_\ell(\omega)}{k} =
\esp{U_0} \biggr\}.
\end{eqnarray*}
By Birkhoff's ergodic theorem, $
\prob{\Omega_2}=1$. To obtain \eqref{eqstationarykplusm}, it is
thus sufficient to show that $\Omega_1^c \cap\Omega_2=\varnothing$. The
proof is by contradiction. Assume $\Omega_1^c \cap\Omega_2 \neq
\varnothing$, so that there exists $\omega\in\Omega_1^c \cap\Omega_2$.
For such $\omega$, the fact that $ \omega\notin\Omega_1$ implies that
there exist a positive number $\epsilon(\omega)>0$ and integer-valued
sequences $(k_n(\omega))_{n \in\nset}$ and $(m_n(\omega))_{n \in\nset
}$ such that $k_n(\omega)+m_n(\omega)\geq n$ and for all $n\geq0$,
%
\begin{equation}
\label{eqcontradic} \biggl\llvert \frac{\sum_{\ell=-m_n(\omega)}^{k_n(\omega)-1}U_{\ell}(\omega
)}{k_n(\omega)+m_n(\omega)}-\esp{U_0} \biggr
\rrvert \geq\epsilon(\omega).
\end{equation}
Consider the following decomposition:
%
\begin{eqnarray}
\label{eqorangis}
\frac{\sum_{\ell=-m_n(\omega)}^{k_n(\omega)-1}U_{\ell}(\omega
)}{k_n(\omega)+m_n(\omega)}
&=& \frac{m_n(\omega)}{k_n(\omega)+m_n(\omega)} \frac{\sum_{\ell
=-m_n(\omega)}^{-1}U_{\ell}(\omega)}{m_n(\omega)}
\nonumber\\[-8pt]\\[-8pt]
&&{} + \frac{k_n(\omega
)}{k_n(\omega)+m_n(\omega)} \frac{\sum_{\ell=0}^{k_n(\omega)-1}U_{\ell
}(\omega)}{k_n(\omega)}.\nonumber
\end{eqnarray}
First, assume that $(k_n(\omega))_{n \in\nset}$ is bounded. Since
$k_n(\omega)+m_n(\omega)\geq n$, it follows that $m_n(\omega)$ tends to
infinity, implying that
%
\begin{eqnarray}\label{eqboisdelepine}
\lim_{n \to\infty}\frac{m_n(\omega)}{k_n(\omega)+m_n(\omega)}&=&1,
\nonumber\\[-8pt]\\[-8pt]
\lim _{n \to\infty}\frac{k_n(\omega)}{k_n(\omega)+m_n(\omega)}&=&0,\nonumber
\end{eqnarray}
whereas $\sum_{\ell=0}^{k_n(\omega)-1}U_{\ell}(\omega)/k_n(\omega)$
remains bounded. However, since $\omega\in\Omega_2$ and $\lim_{n \to
\infty} m_n(\omega) = \infty$,
\[
\lim_{n \to\infty}\frac{\sum_{\ell=-m_n(\omega)}^{-1}U_{\ell}(\omega
)}{m_n(\omega)}=\esp{U_0},
\]
which implies, together with \eqref{eqboisdelepine}, that
\[
\lim_{n \to\infty} \frac{\sum_{\ell=-m_n(\omega)}^{k_n(\omega
)-1}U_{\ell}(\omega)}{k_n(\omega)+m_n(\omega)}= \esp{U_0}.
\]
This contradicts \eqref{eqcontradic}. Using similar arguments one
proves that $(m_n(\omega))_{n \in\nset}$ is unbounded as well. Hence,
we have proved that neither $(k_n(\omega))_{n \in\nset}$ nor
$(m_n(\omega))_{n \in\nset}$ are bounded.

Then, by extracting a subsequence if necessary, one may assume that
$\lim_{n \to\infty} k_n(\omega)=\lim_{n \to\infty} m_n(\omega)=\infty
$. Since $\omega\in\Omega_2$, this implies that
\[
\lim_{n \to\infty} \frac{\sum_{\ell=-m_n(\omega)}^{-1}U_{\ell}(\omega
)}{m_n(\omega)}=\lim_{n \to\infty}
\frac{\sum_{\ell=0}^{k_n(\omega
)-1}U_{\ell}(\omega)}{k_n(\omega)}=\esp{U_0}.  %
\]
Combining this with \eqref{eqorangis}, we obtain that
\[
\lim_{n \to\infty} \frac{\sum_{\ell=-m_n(\omega)}^{k_n(\omega
)-1}U_{\ell}(\omega)}{k_n(\omega)+m_n(\omega)}=\esp{U_0},
\]
which again contradicts \eqref{eqcontradic}. Finally, $\Omega_1^c \cap
\Omega_2=\varnothing$, and since $
\prob{\Omega_2}=1$, we finally obtain that $
\prob{\Omega_1}=1$. The proof is complete.
\end{pf}
%
%
\begin{lemma} \label{lemboundFilter}
Let $\sequence{U}[k][\zset]$, $\sequence{V}[k][\zset]$ and $\sequence
{W}[k][\zset]$ be stationary sequences such that
\[
\espbig{\ln^+ U_0}<\infty,  \qquad\espbig{\ln^+ V_0} < \infty, \qquad\espbig {\ln^+ W_0} < \infty.  %
\]
Then for all $\eta$ and $\rho$ in $\ooint{0,1}$ such that $- \ln\eta>
\espTxt{\ln^+ V_0}$ there exist a $\probalone$\mbox{-a.s.} finite random
variable $C$ and a constant $\beta\in\,\ooint{0,1}$ such that for all
$k \in\nsetpos$ and $m \in\nset$, $\probalone$\mbox{-a.s.},
\[
\rho^{k + m} + \eta^{k + m} W_{-m} \Biggl( \prod
_{\ell= -m}^{k - 1} V_\ell \Biggr)
U_k \leq C \beta^{k + m}.  %
\]
\end{lemma}
\begin{pf}
See \cite{doucmoulines2011}, Lemma 6.
\end{pf}
\end{appendix}

\section*{Acknowledgment}
We thank the anonymous referee for insightful comments that improved
the presentation of the paper.


%

\printaddresses


\begin{thebibliography}{34}

\bibitem{baincrisan2009}
%
\begin{bbook}[mr]
\bauthor{\bsnm{Bain},~\bfnm{Alan}\binits{A.}} \AND
\bauthor{\bsnm{Crisan},~\bfnm{Dan}\binits{D.}}
(\byear{2009}).
\btitle{Fundamentals of Stochastic Filtering}.
\bseries{Stochastic Modelling and Applied Probability}
\bvolume{60}.
\bpublisher{Springer},
\blocation{New York}.
\bid{mr={2454694}}
\end{bbook}
%
\bptok{imsref}%
\endbibitem

\bibitem{cappemoulinesryden2005}
%
\begin{bbook}[mr]
\bauthor{\bsnm{Capp{\'e}},~\bfnm{Olivier}\binits{O.}},
\bauthor{\bsnm{Moulines},~\bfnm{Eric}\binits{E.}} \AND
\bauthor{\bsnm{Ryd{\'e}n},~\bfnm{Tobias}\binits{T.}}
(\byear{2005}).
\btitle{Inference in Hidden {M}arkov Models}.
\bpublisher{Springer},
\blocation{New York}.
\bid{mr={2159833}}
\end{bbook}
%
\bptok{imsref}%
\endbibitem

\bibitem{chopin2002}
%
\begin{barticle}[mr]
\bauthor{\bsnm{Chopin},~\bfnm{Nicolas}\binits{N.}}
(\byear{2002}).
\btitle{A sequential particle filter method for static models}.
\bjournal{Biometrika}
\bvolume{89}
\bpages{539--551}.
\bid{doi={10.1093/biomet/89.3.539}, mr={1929161}}
\end{barticle}
%
\bptok{imsref}%
\endbibitem

\bibitem{crisanheine2008}
%
\begin{barticle}[mr]
\bauthor{\bsnm{Crisan},~\bfnm{D.}\binits{D.}} \AND
\bauthor{\bsnm{Heine},~\bfnm{K.}\binits{K.}}
(\byear{2008}).
\btitle{Stability of the discrete time filter in terms of the tails of
noise distributions}.
\bjournal{J. Lond. Math. Soc. (2)}
\bvolume{78}
\bpages{441--458}.
\bid{doi={10.1112/jlms/jdn032}, mr={2439634}}
\end{barticle}
%
\bptok{imsref}%
\endbibitem

\bibitem{delmoral2004}
%
\begin{bbook}[mr]
\bauthor{\bparticle{Del}~\bsnm{Moral},~\bfnm{Pierre}\binits{P.}}
(\byear{2004}).
\btitle{Feynman--{K}ac Formulae. Genealogical and Interacting Particle Systems with Applications}.
\bpublisher{Springer},
\blocation{New York}.
\bid{doi={10.1007/978-1-4684-9393-1}, mr={2044973}}
\end{bbook}
%
\bptok{imsref}%
\endbibitem

\bibitem{delmoralguionnet1999}
%
\begin{barticle}[mr]
\bauthor{\bparticle{Del}~\bsnm{Moral},~\bfnm{P.}\binits{P.}} \AND
\bauthor{\bsnm{Guionnet},~\bfnm{A.}\binits{A.}}
(\byear{1999}).
\btitle{Central limit theorem for nonlinear filtering and interacting particle systems}.
\bjournal{Ann. Appl. Probab.}
\bvolume{9}
\bpages{275--297}.
\bid{doi={10.1214/aoap/1029962742}, mr={1687359}}
\end{barticle}
%
\bptok{imsref}%
\endbibitem

\bibitem{delmoralguionnet2001}
%
\begin{barticle}[mr]
\bauthor{\bparticle{Del}~\bsnm{Moral},~\bfnm{Pierre}\binits{P.}} \AND
\bauthor{\bsnm{Guionnet},~\bfnm{Alice}\binits{A.}}
(\byear{2001}).
\btitle{On the stability of interacting processes with applications to filtering and genetic algorithms}.
\bjournal{Ann. Inst. Henri Poincar\'e Probab. Stat.}
\bvolume{37}
\bpages{155--194}.
\bid{doi={10.1016/S0246-0203(00)01064-5}, mr={1819122}}
\end{barticle}
%
\bptok{imsref}%
\endbibitem

\bibitem{delmoraljacodprotter2001}
%
\begin{barticle}[mr]
\bauthor{\bparticle{Del} \bsnm{Moral},~\bfnm{Pierre}\binits{P.}},
\bauthor{\bsnm{Jacod},~\bfnm{Jean}\binits{J.}} \AND
\bauthor{\bsnm{Protter},~\bfnm{Philip}\binits{P.}}
(\byear{2001}).
\btitle{The {M}onte-{C}arlo method for filtering with discrete-time
observations}.
\bjournal{Probab. Theory Related Fields}
\bvolume{120}
\bpages{346--368}.
\bid{doi={10.1007/PL00008786}, mr={1843179}}
\end{barticle}
%
\bptok{imsref}%
\endbibitem

\bibitem{delmoralledoux2000}
%
\begin{barticle}[mr]
\bauthor{\bsnm{Del Moral},~\bfnm{P.}\binits{P.}} \AND
\bauthor{\bsnm{Ledoux},~\bfnm{M.}\binits{M.}}
(\byear{2000}).
\btitle{Convergence of empirical processes for interacting particle
systems with applications to nonlinear filtering}.
\bjournal{J. Theoret. Probab.}
\bvolume{13}
\bpages{225--257}.
\bid{doi={10.1023/A:1007743111861}, mr={1744985}}
\end{barticle}
%
\bptok{imsref}%
\endbibitem

\bibitem{delmoralmiclo2000}
%
\begin{bincollection}[mr]
\bauthor{\bparticle{Del}~\bsnm{Moral},~\bfnm{P.}\binits{P.}} \AND
\bauthor{\bsnm{Miclo},~\bfnm{L.}\binits{L.}}
(\byear{2000}).
\btitle{Branching and interacting particle systems approximations of
{F}eynman--{K}ac formulae with applications to nonlinear filtering}.
In \bbooktitle{S\'eminaire de {P}robabilit\'es, {XXXIV}}.
\bseries{Lecture Notes in Math.}
\bvolume{1729}
\bpages{1--145}.
\bpublisher{Springer},
\blocation{Berlin}.
\bid{doi={10.1007/BFb0103798}, mr={1768060}}
\end{bincollection}\vadjust{\goodbreak}
%
\bptok{imsref}%
\endbibitem

\bibitem{doucfortmoulinespriouret2009}
%
\begin{barticle}[mr]
\bauthor{\bsnm{Douc},~\bfnm{R.}\binits{R.}},
\bauthor{\bsnm{Fort},~\bfnm{G.}\binits{G.}},
\bauthor{\bsnm{Moulines},~\bfnm{E.}\binits{E.}} \AND
\bauthor{\bsnm{Priouret},~\bfnm{P.}\binits{P.}}
(\byear{2009}).
\btitle{Forgetting the initial distribution for hidden {M}arkov models}.
\bjournal{Stochastic Process. Appl.}
\bvolume{119}
\bpages{1235--1256}.
\bid{doi={10.1016/j.spa.2008.05.007}, mr={2508572}}
\end{barticle}
%
\bptok{imsref}%
\endbibitem

\bibitem{doucgariviermoulinesolsson2009}
%
\begin{barticle}[mr]
\bauthor{\bsnm{Douc},~\bfnm{Randal}\binits{R.}},
\bauthor{\bsnm{Garivier},~\bfnm{Aur{\'e}lien}\binits{A.}},
\bauthor{\bsnm{Moulines},~\bfnm{Eric}\binits{E.}} \AND
\bauthor{\bsnm{Olsson},~\bfnm{Jimmy}\binits{J.}}
(\byear{2011}).
\btitle{Sequential {M}onte {C}arlo smoothing for general state space
hidden {M}arkov models}.
\bjournal{Ann. Appl. Probab.}
\bvolume{21}
\bpages{2109--2145}.
\bid{doi={10.1214/10-AAP735}, mr={2895411}}
\end{barticle}
%
\bptok{imsref}%
\endbibitem

\bibitem{doucguillinnajim2005}
%
\begin{barticle}[mr]
\bauthor{\bsnm{Douc},~\bfnm{R.}\binits{R.}},
\bauthor{\bsnm{Guillin},~\bfnm{A.}\binits{A.}} \AND
\bauthor{\bsnm{Najim},~\bfnm{J.}\binits{J.}}
(\byear{2005}).
\btitle{Moderate deviations for particle filtering}.
\bjournal{Ann. Appl. Probab.}
\bvolume{15}
\bpages{587--614}.
\bid{doi={10.1214/105051604000000657}, mr={2114983}}
\end{barticle}
%
\bptok{imsref}%
\endbibitem

\bibitem{doucmoulines2008}
%
\begin{barticle}[mr]
\bauthor{\bsnm{Douc},~\bfnm{Randal}\binits{R.}} \AND
\bauthor{\bsnm{Moulines},~\bfnm{Eric}\binits{E.}}
(\byear{2008}).
\btitle{Limit theorems for weighted samples with applications to
sequential {M}onte {C}arlo methods}.
\bjournal{Ann. Statist.}
\bvolume{36}
\bpages{2344--2376}.
\bid{doi={10.1214/07-AOS514}, mr={2458190}}
\end{barticle}
%
\bptok{imsref}%
\endbibitem

\bibitem{doucmoulines2011}
%
\begin{barticle}[mr]
\bauthor{\bsnm{Douc},~\bfnm{Randal}\binits{R.}} \AND
\bauthor{\bsnm{Moulines},~\bfnm{Eric}\binits{E.}}
(\byear{2012}).
\btitle{Asymptotic properties of the maximum likelihood estimation in
misspecified hidden {M}arkov models}.
\bjournal{Ann. Statist.}
\bvolume{40}
\bpages{2697--2732}.
\bid{doi={10.1214/12-AOS1047}, mr={3097617}}
\end{barticle}
%
\bptok{imsref}%
\endbibitem

\bibitem{doucmoulinesolsson2008}
%
\begin{barticle}[mr]
\bauthor{\bsnm{Douc},~\bfnm{Randal}\binits{R.}},
\bauthor{\bsnm{Moulines},~\bfnm{{\'E}ric}\binits{{\'E}.}} \AND
\bauthor{\bsnm{Olsson},~\bfnm{Jimmy}\binits{J.}}
(\byear{2009}).
\btitle{Optimality of the auxiliary particle filter}.
\bjournal{Probab. Math. Statist.}
\bvolume{29}
\bpages{1--28}.
\bid{mr={2552996}}
\end{barticle}
%
\bptok{imsref}%
\endbibitem

\bibitem{doucetdefreitasgordon2001}
%
\begin{bbook}[mr]
\beditor{\bsnm{Doucet} \bfnm{Arnaud}\binits{A.}},
\beditor{\bparticle{de} \bsnm{Freitas} \bfnm{Nando}\binits{N.}} \AND
\beditor{\bsnm{Gordon} \bfnm{Neil}\binits{N.}}, eds.
(\byear{2001}).
\btitle{Sequential {M}onte {C}arlo Methods in Practice}.
\bpublisher{Springer},
\blocation{New York}.
\bid{mr={1847783}}
\end{bbook}
%
\bptok{imsref}%
\endbibitem

\bibitem{godsilldoucetwest2004}
%
\begin{barticle}[auto:STB|2014/02/12|12:18:25]
\bauthor{\bsnm{Godsill},~\bfnm{S.~J.}\binits{S.~J.}},
\bauthor{\bsnm{Doucet},~\bfnm{A.}\binits{A.}} \AND
\bauthor{\bsnm{West},~\bfnm{M.}\binits{M.}}
(\byear{2004}).
\btitle{Monte Carlo smoothing for nonlinear time series}.
\bjournal{J. Amer. Statist. Assoc.}
\bvolume{50}
\bpages{438--449}.
\end{barticle}
%
\bptok{imsref}%
\endbibitem

\bibitem{gordonsalmondsmith1993}
%
\begin{barticle}[auto:STB|2014/02/12|12:18:25]
\bauthor{\bsnm{Gordon},~\bfnm{N.}\binits{N.}},
\bauthor{\bsnm{Salmond},~\bfnm{D.}\binits{D.}} \AND
\bauthor{\bsnm{Smith},~\bfnm{A.~F.}\binits{A.~F.}}
(\byear{1993}).
\btitle{Novel approach to nonlinear/non-{G}aussian Bayesian state estimation}.
\bjournal{IEE Proc., F, Radar Signal Process.}
\bvolume{140}
\bpages{107--113}.
\end{barticle}
%
\bptok{imsref}%
\endbibitem

\bibitem{heinecrisan2008}
%
\begin{barticle}[mr]
\bauthor{\bsnm{Heine},~\bfnm{Kari}\binits{K.}} \AND
\bauthor{\bsnm{Crisan},~\bfnm{Dan}\binits{D.}}
(\byear{2008}).
\btitle{Uniform approximations of discrete-time filters}.
\bjournal{Adv. in Appl. Probab.}
\bvolume{40}
\bpages{979--1001}.
\bid{mr={2488529}}
\end{barticle}
%
\bptok{imsref}%
\endbibitem

\bibitem{doucetjohansen2008}
%
\begin{barticle}[mr]
\bauthor{\bsnm{Johansen},~\bfnm{Adam~M.}\binits{A.~M.}} \AND
\bauthor{\bsnm{Doucet},~\bfnm{Arnaud}\binits{A.}}
(\byear{2008}).
\btitle{A note on auxiliary particle filters}.
\bjournal{Statist. Probab. Lett.}
\bvolume{78}
\bpages{1498--1504}.
\bid{doi={10.1016/j.spl.2008.01.032}, mr={2528342}}
\end{barticle}
%
\bptok{imsref}%
\endbibitem

\bibitem{kingman1973}
%
\begin{barticle}[mr]
\bauthor{\bsnm{Kingman},~\bfnm{J.~F.~C.}\binits{J.~F.~C.}}
(\byear{1973}).
\btitle{Subadditive ergodic theory}.
\bjournal{Ann. Probab.}
\bvolume{1}
\bpages{883--909}.
\bid{mr={0356192}}
\end{barticle}
%
\bptok{imsref}%
\endbibitem

\bibitem{kleptsynaveretennikov2008}
%
\begin{barticle}[mr]
\bauthor{\bsnm{Kleptsyna},~\bfnm{M.~L.}\binits{M.~L.}} \AND
\bauthor{\bsnm{Veretennikov},~\bfnm{A.~Yu.}\binits{A.~Yu.}}
(\byear{2008}).
\btitle{On discrete time ergodic filters with wrong initial data}.
\bjournal{Probab. Theory Related Fields}
\bvolume{141}
\bpages{411--444}.
\bid{doi={10.1007/s00440-007-0089-7}, mr={2391160}}
\end{barticle}
%
\bptok{imsref}%
\endbibitem

\bibitem{kuensch2005}
%
\begin{barticle}[mr]
\bauthor{\bsnm{K{\"u}nsch},~\bfnm{Hans~R.}\binits{H.~R.}}
(\byear{2005}).
\btitle{Recursive {M}onte {C}arlo filters: Algorithms and theoretical analysis}.
\bjournal{Ann. Statist.}
\bvolume{33}
\bpages{1983--2021}.
\bid{doi={10.1214/009053605000000426}, mr={2211077}}
\end{barticle}
%
\bptok{imsref}%
\endbibitem

\bibitem{leglandoudjane2003}
%
\begin{barticle}[mr]
\bauthor{\bsnm{LeGland},~\bfnm{Fran{\c{c}}ois}\binits{F.}} \AND
\bauthor{\bsnm{Oudjane},~\bfnm{Nadia}\binits{N.}}
(\byear{2003}).
\btitle{A robustification approach to stability and to uniform particle
approximation of nonlinear filters: The example of pseudo-mixing signals}.
\bjournal{Stochastic Process. Appl.}
\bvolume{106}
\bpages{279--316}.
\bid{doi={10.1016/S0304-4149(03)00041-3}, mr={1989630}}
\end{barticle}
%
\bptok{imsref}%
\endbibitem

\bibitem{leglandoudjane2004}
%
\begin{barticle}[mr]
\bauthor{\bsnm{Le Gland},~\bfnm{Fran{\c{c}}ois}\binits{F.}} \AND
\bauthor{\bsnm{Oudjane},~\bfnm{Nadia}\binits{N.}}
(\byear{2004}).
\btitle{Stability and uniform approximation of nonlinear filters using
the {H}ilbert metric and application to particle filters}.
\bjournal{Ann. Appl. Probab.}
\bvolume{14}
\bpages{144--187}.
\bid{doi={10.1214/aoap/1075828050}, mr={2023019}}
\end{barticle}
%
\bptok{imsref}%
\endbibitem

\bibitem{olssonryden2008}
%
\begin{barticle}[mr]
\bauthor{\bsnm{Olsson},~\bfnm{Jimmy}\binits{J.}} \AND
\bauthor{\bsnm{Ryd{\'e}n},~\bfnm{Tobias}\binits{T.}}
(\byear{2008}).
\btitle{Asymptotic properties of particle filter-based maximum
likelihood estimators for state space models}.
\bjournal{Stochastic Process. Appl.}
\bvolume{118}
\bpages{649--680}.
\bid{doi={10.1016/j.spa.2007.05.007}, mr={2394847}}
\end{barticle}
%
\bptok{imsref}%
\endbibitem

\bibitem{oudjanerubenthaler2005}
%
\begin{barticle}[mr]
\bauthor{\bsnm{Oudjane},~\bfnm{Nadia}\binits{N.}} \AND
\bauthor{\bsnm{Rubenthaler},~\bfnm{Sylvain}\binits{S.}}
(\byear{2005}).
\btitle{Stability and uniform particle approximation of nonlinear
filters in case of nonergodic signals}.
\bjournal{Stoch. Anal. Appl.}
\bvolume{23}
\bpages{421--448}.
\bid{doi={10.1081/SAP-200056643}, mr={2140972}}
\end{barticle}
%
\bptok{imsref}%
\endbibitem

\bibitem{pittshephard1999}
%
\begin{barticle}[mr]
\bauthor{\bsnm{Pitt},~\bfnm{Michael~K.}\binits{M.~K.}} \AND
\bauthor{\bsnm{Shephard},~\bfnm{Neil}\binits{N.}}
(\byear{1999}).
\btitle{Filtering via simulation: Auxiliary particle filters}.
\bjournal{J. Amer. Statist. Assoc.}
\bvolume{94}
\bpages{590--599}.
\bid{doi={10.2307/2670179}, mr={1702328}}
\end{barticle}
%
\bptok{imsref}%
\endbibitem

\bibitem{polya1976a}
%
\begin{bbook}[mr]
\bauthor{\bsnm{P{\'o}lya},~\bfnm{G.}\binits{G.}} \AND
\bauthor{\bsnm{Szeg{\H{o}}},~\bfnm{G.}\binits{G.}}
(\byear{1976}).
\btitle{Problems and Theorems in Analysis}.
\bpublisher{Springer},
\blocation{New York}.
\bid{mr={0396134}}
\end{bbook}
%
\bptok{imsref}%
\endbibitem

\bibitem{serfling1980}
%
\begin{bbook}[mr]
\bauthor{\bsnm{Serfling},~\bfnm{Robert~J.}\binits{R.~J.}}
(\byear{1980}).
\btitle{Approximation Theorems of Mathematical Statistics}.
\bpublisher{Wiley},
\blocation{New York}.
\bid{mr={0595165}}
\end{bbook}
%
\bptok{imsref}%
\endbibitem

\bibitem{shiryaev1996}
%
\begin{bbook}[mr]
\bauthor{\bsnm{Shiryaev},~\bfnm{A.~N.}\binits{A.~N.}}
(\byear{1996}).
\btitle{Probability},
\bedition{2nd} ed.
\bpublisher{Springer},
\blocation{New York}.
\bid{mr={1368405}}
\end{bbook}
%
\bptok{imsref}%
\endbibitem

\bibitem{tadicdoucet2005}
%
\begin{barticle}[mr]
\bauthor{\bsnm{Tadi{\'c}},~\bfnm{Vladislav~B.}\binits{V.~B.}} \AND
\bauthor{\bsnm{Doucet},~\bfnm{Arnaud}\binits{A.}}
(\byear{2005}).
\btitle{Exponential forgetting and geometric ergodicity for optimal
filtering in general state--space models}.
\bjournal{Stochastic Process. Appl.}
\bvolume{115}
\bpages{1408--1436}.
\bid{doi={10.1016/j.spa.2005.03.005}, mr={2152381}}
\end{barticle}
%
\bptok{imsref}%
\endbibitem

\bibitem{vanhandel2009}
%
\begin{barticle}[mr]
\bauthor{\bparticle{van} \bsnm{Handel},~\bfnm{Ramon}\binits{R.}}
(\byear{2009}).
\btitle{The stability of conditional {M}arkov processes and {M}arkov
chains in random environments}.
\bjournal{Ann. Probab.}
\bvolume{37}
\bpages{1876--1925}.
\bid{doi={10.1214/08-AOP448}, mr={2561436}}
\end{barticle}
%
\bptok{imsref}%
\endbibitem

\end{thebibliography}
\end{document}